# Empirical spectral processes for locally stationary time series

RAINER DAHLHAUS[1] and WOLFGANG POLONIK[2]

[1] *Institut für Angewandte Mathematik, Universität Heidelberg, Im Neuenheimer Feld 294, 69120 Heidelberg, Germany. E-mail: dahlhaus@statlab.uni-heidelberg.de*
[2] *Department of Statistics, University of California, Davis, CA 95616-8705, USA.*
*E-mail: polonik@wald.ucdavis.edu*

A time-varying empirical spectral process indexed by classes of functions is defined for locally stationary time series. We derive weak convergence in a function space, and prove a maximal exponential inequality and a Glivenko–Cantelli-type convergence result. The results use conditions based on the metric entropy of the index class. In contrast to related earlier work, no Gaussian assumption is made. As applications, quasi-likelihood estimation, goodness-of-fit testing and inference under model misspecification are discussed. In an extended application, uniform rates of convergence are derived for local Whittle estimates of the parameter curves of locally stationary time series models.

## 1. Introduction

In recent years, several methods have been derived for locally stationary time series models, that is, for models which can locally be approximated by stationary processes. Out of the large literature, we mention the work of Priestley (1965) on oscillatory processes, Dahlhaus (1997) on locally stationary processes, Neumann and von Sachs (1997) on wavelet estimation of evolutionary spectra, Nason, von Sachs and Kroisandt (2000) on a wavelet-based model of evolutionary spectra and more recent work such as Davis, Lee and Rodriguez-Yam (2005) on piecewise stationary processes, Fryzlewicz, Sapatinas and Subba Rao (2006) on locally stationary volatility estimation and Sakiyama and Taniguchi (2004) on discriminant analysis for locally stationary processes.

In this paper, we emphasize the relevance of the empirical spectral process for locally stationary time series. During the last decade, the theory of empirical processes has developed considerably and the number of statistical problems approached by utilizing concepts from empirical process theory is steeply increasing. In this paper, we show how large parts of the existing methodology on empirical processes can fruitfully be used for time series analysis of locally stationary processes. In our set-up, the role of the empirical distribution of i.i.d. data is taken over by the empirical time-varying spectral measure. This generalizes a similar approach for stationary time series (cf. Dahlhaus (1988), Mikosch and Norvaisa (1997), Fay and Soulier (2001)). An overview of these methods and some references to the existing literature on empirical process techniques in other settings may be found in Dahlhaus and Polonik (2002).

In Section 2, we introduce the empirical spectral process indexed by classes of functions, derive its convergence (including a functional central limit theorem) and prove a maximal exponential inequality and a Glivenko–Cantelli-type convergence result. These results use conditions based on the metric entropy of the index class.

The empirical spectral process plays a key role in many statistical applications. In Section 3, we briefly discuss parametric quasi-likelihood estimation, nonparametric quasi-likelihood estimation, inference under model misspecification by stationary models and local estimates. An extended application is given in Section 4, where uniform rates of convergence are derived for local Whittle estimates of the parameter curves of locally stationary time series models.

Although our concept is based on empirical process techniques in the frequency domain, there exist many applications in the time domain. Section 3 and Section 4 contain many examples in the time domain, particularly with time-varying ARMA models.

The empirical spectral process for locally stationary processes has also been briefly considered in Dahlhaus and Polonik (2006) in the special context of nonparametric estimation. In comparison to that paper, we also consider here the case of non-Gaussian processes and use weaker assumptions on the underlying process. We mention that the assumptions on the underlying process are very weak, allowing for jumps in the parameter curves by assuming bounded variation instead of continuity in the time direction.

All proofs are deferred without further reference to Section 5.

## 2. The time-varying empirical spectral process

In this section, we define the empirical spectral process and derive its properties including a functional central limit theorem and a maximal exponential inequality.

### 2.1. Locally stationary processes

Locally stationary processes were introduced in Dahlhaus (1997) by using a time-varying spectral representation. In contrast to this, in this paper, we use a time-varying $MA(\infty)$-representation and formulate the assumptions in the time domain. As in nonparametric regression, we rescale the functions in time to the unit interval in order to achieve a meaningful asymptotic theory. The following assumptions on the locally stationary process are the same as those used in Dahlhaus and Polonik (2006). They are more general than, for example, in Dahlhaus (1997) since the parameter curves are allowed to have jumps.

Let

$$V(g) = \sup\left\{\sum_{k=1}^{m} |g(x_k) - g(x_{k-1})| : 0 \le x_0 < \cdots < x_m \le 1,\, m \in \mathbf{N}\right\} \tag{1}$$

be the total variation of a function $g$ on $[0,1]$ and for some $\kappa > 0$, let

$$\ell(j) := \begin{cases} 1, & |j| \le 1, \\ |j| \log^{1+\kappa} |j|, & |j| > 1. \end{cases}$$

**Assumption 2.1.** $X_{t,n}$ $(t = 1, \ldots, n)$ *has a representation*

$$X_{t,n} = \sum_{j=-\infty}^{\infty} a_{t,n}(j)\,\varepsilon_{t-j} \tag{2}$$

*satisfying the following conditions:*

$$\sup_{t,n} |a_{t,n}(j)| \le \frac{K}{\ell(j)} \tag{3}$$

*and there exist functions* $a(\cdot, j) : (0, 1] \to \mathbf{R}$ *with*

$$\sup_u |a(u, j)| \le \frac{K}{\ell(j)}, \tag{4}$$

$$\sup_j \sum_{t=1}^{n} \left| a_{t,n}(j) - a\left(\frac{t}{n}, j\right) \right| \le K, \tag{5}$$

$$V(a(\cdot, j)) \le \frac{K}{\ell(j)}. \tag{6}$$

*The* $\varepsilon_t$ *are assumed to be independent and identically distributed with* $\mathbf{E}\varepsilon_t \equiv 0$ *and* $\mathbf{E}\varepsilon_t^2 \equiv 1$. *In addition, we assume that all moments of* $\varepsilon_t$ *exist and set* $\kappa_4 := \mathrm{cum}_4(\varepsilon_t)$.

**Remark 2.2.** (i) The rather complicated construction with different coefficients $a_{t,n}(j)$ and $a(\frac{t}{n}, j)$ is necessary since we need, on the one hand, a certain smoothness in the time direction (guaranteed by bounded variation of the functions $a(u, j)$) and on the other hand, a class which is rich enough to cover interesting examples. For instance, Proposition 2.4 implies that the process $X_{t,n} = \phi(\frac{t}{n}) X_{t-1,n} + \varepsilon_t$ has a representation of the form (1). However, the proof of Proposition 2.4 reveals that this $X_{t,n}$ does not have a representation of the form

$$X_{t,n} = \sum_{j=-\infty}^{\infty} a\left(\frac{t}{n}, j\right) \varepsilon_{t-j}. \tag{7}$$

(ii) The time-varying MA($\infty$)-representation (2) can easily be transformed into a time-varying spectral representation as used, for example, in Dahlhaus (1997). If the $\varepsilon_t$ are assumed to be stationary, then there exists a Cramér representation

$$\varepsilon_t = \frac{1}{\sqrt{2\pi}} \int_{-\pi}^{\pi} \exp(\mathrm{i}\lambda t)\, \mathrm{d}\xi(\lambda),$$

where $\xi(\lambda)$ is a process with mean 0 and orthonormal increments (cf. Brillinger (1981)). Let

$$A_{t,n}(\lambda) := \sum_{j=-\infty}^{\infty} a_{t,n}(j) \exp(-\mathrm{i}\lambda j). \tag{8}$$

Then

$$X_{t,n} = \frac{1}{\sqrt{2\pi}} \int_{-\pi}^{\pi} \exp(\mathrm{i}\lambda t) A_{t,n}(\lambda)\, \mathrm{d}\xi(\lambda). \tag{9}$$

If (5) is replaced by the stronger condition

$$\sup_t \left| a_{t,n}(j) - a\left(\frac{t}{n}, j\right) \right| \leq \frac{K}{n\ell(j)},$$

then it follows that

$$\sup_{t,\lambda} \left| A_{t,n}(\lambda) - A\left(\frac{t}{n}, \lambda\right) \right| \leq K n^{-1}, \tag{10}$$

which was assumed in Dahlhaus (1997). Conversely, if we start with (9) and (10), then the conditions of Assumption 2.1 can be derived from adequate smoothness conditions on $A(u, \lambda)$.

**Definition 2.3 (Time-varying spectral density and covariance).** *The function*

$$f(u, \lambda) := \frac{1}{2\pi} |A(u, \lambda)|^2$$

*with*

$$A(u, \lambda) := \sum_{j=-\infty}^{\infty} a(u, j) \exp(-\mathrm{i}\lambda j)$$

*is the time-varying spectral density and*

$$c(u, k) := \int_{-\pi}^{\pi} f(u, \lambda) \exp(\mathrm{i}\lambda k)\, \mathrm{d}\lambda = \sum_{j=-\infty}^{\infty} a(u, k+j) a(u, j) \tag{11}$$

*is the time-varying covariance of lag $k$ at rescaled time $u$.*

For a deeper understanding of the time-varying covariance, see also Proposition 5.4.

A simple example of a process $X_{t,n}$ which fulfills the above assumptions is $X_{t,n} = \phi(\frac{t}{n})Y_t$, where $Y_t = \sum_j a(j)\varepsilon_{t-j}$ is stationary with $|a(j)| \leq K/\ell(j)$ and $\phi$ is of bounded variation. From the following proposition, it follows that time-varying ARMA (tvARMA) models whose coefficient functions are of bounded variation are locally stationary in the above sense. The result is proved in Appendix.

**Proposition 2.4 (tvARMA).** *Consider the system of difference equations*

$$\sum_{j=0}^{p} \alpha_j\left(\frac{t}{n}\right) X_{t-j,n} = \sum_{k=0}^{q} \beta_k\left(\frac{t}{n}\right)\sigma\left(\frac{t-k}{n}\right)\varepsilon_{t-k}, \tag{12}$$

*where $\varepsilon_t$ are i.i.d. with $\mathbf{E}\varepsilon_t = 0$, $\mathbf{E}|\varepsilon_t| < \infty$, $\alpha_0(u) \equiv \beta_0(u) \equiv 1$ and $\alpha_j(u) = \alpha_j(0)$, $\beta_k(u) = \beta_k(0)$ for $u < 0$. If all $\alpha_j(\cdot)$ and $\beta_k(\cdot)$, as well as $\sigma^2(\cdot)$, are of bounded variation and $\sum_{j=0}^{p}\alpha_j(u)z^j \neq 0$ for all $u$ and all $0 < |z| \leq 1 + \delta$ for some $\delta > 0$, then there exists a solution of the form*

$$X_{t,n} = \sum_{j=0}^{\infty} a_{t,n}(j)\,\varepsilon_{t-j}$$

*which fulfills (3)–(6) of Assumption 2.1. The time-varying spectral density is given by*

$$f(u,\lambda) = \frac{\sigma^2(u)}{2\pi} \frac{|\sum_{k=0}^{q}\beta_k(u)\exp(\mathrm{i}\lambda k)|^2}{|\sum_{j=0}^{p}\alpha_j(u)\exp(\mathrm{i}\lambda j)|^2}.$$

## 2.2. Convergence of the empirical spectral process

The empirical spectral process is defined by

$$E_n(\phi) = \sqrt{n}(F_n(\phi) - F(\phi)), \tag{13}$$

where

$$F(\phi) = \int_0^1 \int_{-\pi}^{\pi} \phi(u,\lambda) f(u,\lambda)\,\mathrm{d}\lambda\,\mathrm{d}u \tag{14}$$

and

$$F_n(\phi) = \frac{1}{n}\sum_{t=1}^{n}\int_{-\pi}^{\pi}\phi\left(\frac{t}{n},\lambda\right)J_n\left(\frac{t}{n},\lambda\right)\mathrm{d}\lambda \tag{15}$$

with the pre-periodogram

$$J_n\left(\frac{t}{n},\lambda\right) = \frac{1}{2\pi}\sum_{k:1\leq[t+1/2\pm k/2]\leq n} X_{[t+1/2+k/2],n}X_{[t+1/2-k/2],n}\exp(-\mathrm{i}\lambda k). \tag{16}$$

If $X_{[t+1/2+k/2],n}X_{[t+1/2-k/2],n}$ is regarded as a (raw) estimate of $c(\frac{t}{n},k)$, then $J_n(\frac{t}{n},\lambda)$ can be regarded as a (raw) estimate of $f(\frac{t}{n},\lambda)$. However, in order to become consistent, $J_n(\frac{t}{n},\lambda)$ needs to be smoothed in time and frequency directions. The pre-periodogram $J_n$ was first defined by Neumann and von Sachs (1997).

Many statistics occurring in the analysis of non-stationary time series can be written as functionals of $F_n(\phi)$. Several examples are discussed in Section 3 and Section 4.

We first prove a central limit theorem for $E_n(\phi)$ under the assumption that we have bounded variation in both components of $\phi(u,\lambda)$. Besides the definition in (1), we need a definition in two dimensions. Let

$$V^2(\phi) = \sup\left\{ \sum_{j,k=1}^{\ell,m} |\phi(u_j,\lambda_k) - \phi(u_{j-1},\lambda_k) - \phi(u_j,\lambda_{k-1}) + \phi(u_{j-1},\lambda_{k-1})| : \right.$$

$$\left. 0 \le u_0 < \cdots < u_\ell \le 1; -\pi \le \lambda_0 < \cdots < \lambda_m \le \pi; \ell, m \in \mathbf{N} \right\}.$$

For simplicity, we set

$$\|\phi\|_{\infty,V} := \sup_u V(\phi(u,\cdot)), \qquad \|\phi\|_{V,\infty} := \sup_\lambda V(\phi(\cdot,\lambda)),$$

$$\|\phi\|_{V,V} := V^2(\phi) \quad \text{and} \quad \|\phi\|_{\infty,\infty} := \sup_{u,\lambda} |\phi(u,\lambda)|.$$

**Theorem 2.5.** *Suppose Assumption 2.1 is fulfilled and $\phi_1,\ldots,\phi_k$ are functions with $\|\phi_j\|_{\infty,V}$, $\|\phi_j\|_{V,\infty}$, $\|\phi_j\|_{V,V}$ and $\|\phi_j\|_{\infty,\infty}$ being finite $(j=1,\ldots,k)$. Then*

$$(E_n(\phi_j))_{j=1,\ldots,k} \xrightarrow{\mathcal{D}} (E(\phi_j))_{j=1,\ldots,k},$$

*where $(E(\phi_j))_{j=1,\ldots,k}$ is a Gaussian random vector with mean 0 and*

$$\mathrm{cov}(E(\phi_j),E(\phi_k))$$
$$= 2\pi \int_0^1 \int_{-\pi}^\pi \phi_j(u,\lambda)\left[\phi_k(u,\lambda) + \phi_k(u,-\lambda)\right] f^2(u,\lambda)\,\mathrm{d}\lambda\,\mathrm{d}u$$
$$+ \kappa_4 \int_0^1 \left(\int_{-\pi}^\pi \phi_j(u,\lambda_1)f(u,\lambda_1)\,\mathrm{d}\lambda_1\right)\left(\int_{-\pi}^\pi \phi_k(u,\lambda_2)f(u,\lambda_2)\,\mathrm{d}\lambda_2\right)\mathrm{d}u.$$

**Remark 2.6.** (i) We mention that Theorem 5.3 contains a similar statement under a different set of conditions which is obtained as a by-product of our calculations. Furthermore, we mention that Theorem 2.5 also holds if a data taper is used, that is, if $F_n(\phi)$ and $F(\phi)$ are defined as in (43) and (42) (in that case, we also need Assumption 5.1 and $\mathrm{cov}(E(\phi_j),E(\phi_k))$ must be replaced by $c_E^{(h)}(\phi_j,\phi_k)$ as defined in Theorem 5.3). For simplicity, we consider the tapered case only in Section 5.

(ii) In contrast to earlier results (cf. Dahlhaus and Neumann (2001), Lemma 2.1), the assumptions on $\phi(u, \lambda)$ and $f(u, \lambda)$ are very weak. In particular, we allow for non-continuous behavior.

(iii) In the stationary case where $\phi_j(u, \lambda) = \tilde{\phi}_j(\lambda)$ and $f(u, \lambda) = \tilde{f}(\lambda)$, this is the classical central limit theorem for the weighted periodogram (see Example 3.3 below).

(iv) The limit behavior for complex-valued $\phi_j$ can easily be derived from Theorem 2.5 by considering the real and imaginary parts separately.

In Theorem 2.11, a functional central limit theorem indexed by function spaces and in Theorem 2.12 a Glivenko–Cantelli-type theorem are proved. The central ingredient of their proofs will be an exponential inequality for the empirical spectral process and a maximal inequality derived in the next subsection.

## 2.3. A maximal exponential inequality

Let

$$
\begin{aligned}
\rho_2(\phi) &:= \left( \int_0^1 \int_{-\pi}^{\pi} \phi(u, \lambda)^2 \, d\lambda \, du \right)^{1/2}, \\
\rho_{2,n}(\phi) &:= \left( \frac{1}{n} \sum_{t=1}^{n} \int_{-\pi}^{\pi} \phi\left( \frac{t}{n}, \lambda \right)^2 d\lambda \right)^{1/2}
\end{aligned}
\tag{17}
$$

and

$$
\tilde{E}_n(\phi) := \sqrt{n} \left( F_n(\phi) - \mathbf{E} F_n(\phi) \right).
\tag{18}
$$

**Theorem 2.7 (Exponential inequality).** *Suppose Assumption 2.1 is fulfilled with $\mathbf{E}|\varepsilon_t|^k \leq C_\varepsilon^k$ for all $k \in \mathbf{N}$. We then have, for all $\eta > 0$,*

$$
P(|\tilde{E}_n(\phi)| \geq \eta) \leq c_1 \exp\left( -c_2 \sqrt{\frac{\eta}{\rho_{2,n}(\phi)}} \right)
\tag{19}
$$

*with some constants $c_1$, $c_2 > 0$ independent of $n$.*

**Remark 2.8.** (i) In the Gaussian case, it is possible to omit the $\sqrt{\cdot}$ in (19) or to prove a Bernstein-type inequality which is even stronger (cf. Dahlhaus and Polonik (2006), Theorem 4.1).

(ii) The assumption $\mathbf{E}|\varepsilon_t|^k \leq C_\varepsilon^k$ for all $k \in \mathbf{N}$ is strong in that it implies finite exponential moments for $\varepsilon_t$. On the other hand, with the above exponential inequality, this leads to a very strong finite-sample result for a large class of locally stationary processes (remember that the assumptions made on the smoothness of the parameter curves are very weak). The strong assumptions on $\varepsilon_t$ make sense for obtaining the uniform bounds of the empirical spectral process and the functional central limit theorem below. Furthermore, they will lead to a very strong uniform rate of convergence result in Theorem 4.1.

(iii) To treat the bias $\mathbf{E}F_n(\phi) - F(\phi)$, we set $F^+(\phi) := \frac{1}{n}\sum_{t=1}^{n}\int_{-\pi}^{\pi}\phi(\frac{t}{n},\lambda)\,f(\frac{t}{n},\lambda)\,\mathrm{d}\lambda$. We then we have (see proof of Theorem 2.7 and Remark 2.8)

$$\sqrt{n}|\mathbf{E}F_n(\phi) - F^+(\phi)| \leq K\rho_{2,n}(\phi), \tag{20}$$

$$\sqrt{n}|F^+(\phi) - F(\phi)| \leq \frac{K}{\sqrt{n}}\left(\|\phi\|_{V,\infty} + \|\phi\|_{\infty,\infty}\right) \tag{21}$$

and

$$\rho_{2,n}(\phi)^2 \leq \rho_2(\phi)^2 + \frac{4\pi}{n}\|\phi\|_{V,\infty}\|\phi\|_{\infty,\infty}, \tag{22}$$

leading, for example, to the exponential inequality

$$P(|E_n(\phi)| \geq \eta) \leq c_1' \exp\left\{-c_2'\frac{\eta^{1/2}}{(\rho_2(\phi) + (1/\sqrt{n})\|\phi\|_{V,\infty} + (1/\sqrt{n})\|\phi\|_{\infty,\infty})^{1/2}}\right\}. \tag{23}$$

An alternative inequality used later is

$$\sqrt{n}|\mathbf{E}F_n(\phi) - F^+(\phi)| \leq K\,\frac{\log n}{\sqrt{n}}\,\|\phi\|_{\infty,V} + \frac{K}{\sqrt{n}}\|\phi\|_{\infty,\infty}. \tag{24}$$

The above exponential inequality is the core of the proof of the following result which then leads to stochastic equicontinuity of the empirical spectral process. Analogously to standard empirical process theory, stochastic equicontinuity is crucial for proving tightness.

As for the standard empirical process, the results for the function-indexed empirical spectral process $(E_n(\phi), \phi \in \Phi)$ are derived under conditions on the richness of $\Phi$, measured by the metric entropy. For each $\epsilon > 0$, the *covering number* of $\Phi$ with respect to the norm $\rho_2$ is defined by

$$N(\epsilon, \Phi, \rho_2) = \inf\{n \geq 1 : \exists\,\phi_1,\ldots,\phi_n \in \Phi \text{ such that}$$

$$\forall\phi \in \Phi\,\exists\,1 \leq i \leq n \text{ with } \rho_2(\phi - \phi_i) \leq \epsilon\}$$

and the *metric entropy* of $\Phi$ with respect to $\rho_2$ By

$$H(\epsilon, \Phi, \rho_2) = \log N(\epsilon, \Phi, \rho_2). \tag{25}$$

Usually, the metric entropy is not known exactly, only upper bounds are known. These upper bounds are usually of the form $H(\epsilon, \Phi, \rho_2) \leq C\epsilon^{-r}$ or $N(\epsilon, \Phi, \rho_2) \leq C\epsilon^{-r}$ with $C, r > 0$. For the results below, we assume

$$H(\epsilon, \Phi, \rho_2) \leq \tilde{H}_\Phi(\epsilon)$$

with an upper bound $\tilde{H}_\Phi(\cdot)$ which is assumed to be continuous and strictly decreasing.

**Remark.** In standard empirical process theory, so-called *bracketing covering numbers* are often used instead. Here, we do not use bracketing covering numbers since the empirical spectral process is not monotone in $\phi$.

Let

$$\tau_{\infty,V} := \sup_{\phi \in \Phi} \|\phi\|_{\infty,V}, \qquad \tau_{V,\infty} := \sup_{\phi \in \Phi} \|\phi\|_{V,\infty},$$

$$\tau_{V,V} := \sup_{\phi \in \Phi} \|\phi\|_{V,V} \quad \text{and} \quad \tau_{\infty,\infty} := \sup_{\phi \in \Phi} \|\phi\|_{\infty,\infty}.$$

In order to avoid further technical assumptions, the following results assume measurability of all random quantities without further mentioning it.

**Theorem 2.9 (Maximal inequality).** *Suppose Assumption 2.1 is fulfilled with* $\mathbf{E}|\varepsilon_t|^k \leq C_\varepsilon^k$ *for all* $k \in \mathbf{N}$. *Suppose that* $\Phi$ *is such that* $\tau_{\infty,V}$, $\tau_{V,\infty}$, $\tau_{V,V}$ *and* $\tau_{\infty,\infty}$ *are finite. Let*

$$\tau_2 := \sup_{\phi \in \Phi} \rho_2(\phi).$$

*There then exists a set* $B_n$ *(independent of* $\Phi$*) with* $\lim_{n \to \infty} P(B_n) = 1$ *and a constant* $L$ *(independent of* $\Phi$ *and* $n$*) such that for all* $\eta$ *satisfying*

$$\eta \geq 26L \max\{\tau_{\infty,V}, \tau_{V,\infty}, \tau_{V,V}, \tau_{\infty,\infty}\} \frac{(\log n)^3}{\sqrt{n}} \tag{26}$$

*and*

$$\eta \geq \frac{72}{c_2^2} \int_0^\alpha \tilde{H}_\Phi(s)^2 \, \mathrm{d}s \qquad \text{with } \alpha := \tilde{H}_\Phi^{-1}\left(\frac{c_2}{4} \sqrt{\frac{\eta}{\tau_2}}\right), \tag{27}$$

*we have*

$$P\left(\sup_{\phi \in \Phi} |\tilde{E}_n(\phi)| > \eta, B_n\right) \leq 3c_1 \exp\left\{-\frac{c_2}{4} \sqrt{\frac{\eta}{\tau_2}}\right\} \tag{28}$$

*and*

$$P\left(\sup_{\phi \in \Phi} |E_n(\phi)| > \eta, B_n\right) \leq 3c_1 \exp\left\{-\frac{c_2}{4} \sqrt{\frac{\eta}{\tau_2}}\right\} \tag{29}$$

*with some constants* $c_1$, $c_2 > 0$ *independent of* $n$.

**Remark 2.10.** (i) A maximal inequality for $\{\tilde{E}_n(\phi), \phi \in \Phi\}$, assuming Gaussian innovations, can be found in Dahlhaus and Polonik (2006). The additional Gaussian assumption enables a weakening of the crucial assumption (27), essentially replacing $\int_0^\alpha \tilde{H}_\Phi^2(s) \, \mathrm{d}s$ by $\int_0^\alpha \tilde{H}_\Phi(s) \, \mathrm{d}s$. Related to that, the resulting exponential inequality is stronger, replacing $\sqrt{\frac{\eta}{\tau}}$ in the exponent of (28) by $\frac{\eta}{\tau}$. The proof in the present non-Gaussian case is much more complicated. It is an open question as to whether the same result as in the Gaussian case can also be obtained for non-Gaussian processes.

(ii) The restriction to the set $B_n$ has several advantages. First, it allows for replacing $\rho_{2,n}(\phi)$ by $\rho_2(\phi)$ (more precisely, by $\tau_2 = \sup_{\phi \in \Phi} \rho_2(\phi)$) which makes the results much simpler. Furthermore, extra terms due to the bias, as in (23), can be avoided. For many results, the set $B_n$ means no restriction, particularly if the probability of an event is calculated as for equicontinuity. The set $B_n$ is given in (72) and it is shown in Lemma 5.9 that $\mathbf{P}(B_n^c) = O(n^{-1})$. For this reason, $B_n$ may be removed from (28) and (29) by adding an $O(n^{-1})$ term to the right-hand side of these quantities. However, for many applications, this would not be sufficient.

(iii) $c_1$ and $c_2$ are the constants from (19). The minimal choice of $L$ is $L = \max\{K_1, K_2, K\}$, where $K_1, K_2 > 0$ are from Lemma 5.8 and $K$ is the constant from (21) and (24).

## 2.4. A functional central limit theorem and a GC-type result

Theorems 2.5 and 2.9 are the main ingredients for deriving the following weak convergence result for the process $\{E_n(\phi); \phi \in \Phi\}$ in the space $\ell^\infty(\Phi)$ of uniformly bounded (real-valued) functions on $\Phi$, that is, with $\|g\|_\Phi := \sup_{\phi \in \Phi} |g(\phi)|$, we have $\ell^\infty(\Phi) = \{g : \Phi \to \mathbf{R}; \|g\|_\Phi < \infty\}$.

**Theorem 2.11 (Functional limit theorem).** *Suppose Assumption 2.1 is fulfilled with* $\mathbf{E}|\varepsilon_t|^k \le C_\varepsilon^k$ *for all* $k \in \mathbf{N}$. *Furthermore let* $\Phi$ *be such that* $\tau_{\infty,V}$, $\tau_{V,\infty}$, $\tau_{V,V}$ *and* $\tau_{\infty,\infty}$ *are finite. If, in addition,*

$$\int_0^1 \tilde{H}_\Phi(s)^2 \, \mathrm{d}s < \infty, \tag{30}$$

*then we have*

$$E_n(\cdot) \to E(\cdot) \qquad weakly\ in\ \ell_\infty(\Phi)$$

*as* $n \to \infty$, *where* $\{E(\phi), \phi \in \Phi\}$ *denotes a tight, mean zero Gaussian process with covariance structure as given in Theorem 2.5.*

Weak convergence in the above theorem means that $\mathbf{E}^* \alpha(E_n) \to \mathbf{E}\alpha(E)$ as $n \to \infty$ for every bounded, continuous real-valued function $\alpha$ on $\ell^\infty(\Phi)$ equipped with the supremum norm, where $\mathbf{E}^*$ denotes outer expectation. This Hoffman–Jørgensen-type formulation of weak convergence avoids measurability considerations for the process $\{E_n(\phi), \phi \in \Phi\}$. Measurability of $E_n$ might become problematic, particularly if $\Phi$ is not separable. Nevertheless, this notion of weak convergence allows the application of useful probabilistic tools such as continuous mapping theorems. For more details, we refer to van der Vaart and Wellner (1996).

Finally, we mention our conjecture that it should be possible to prove another version of the above central limit theorem under much weaker moment assumptions on the $\varepsilon_t$ if the class $\Phi$ is smaller (by avoiding the use of the maximal inequality for proving equicontinuity).

As another application of the maximal inequality, we now prove a Glivenko–Cantelli-type theorem for the empirical spectral process. Here, we allow for a class $\Phi = \Phi_n$ which may be increasing with $n$. We set $\tau_{\infty,V}^{(n)} := \sup_{\phi \in \Phi_n} \|\phi\|_{\infty,V}$, etc.

**Theorem 2.12.** *Suppose Assumption 2.1 is fulfilled with* $\mathbf{E}|\varepsilon_t|^k \leq C_\varepsilon^k$ *for all* $k \in \mathbf{N}$. *Suppose further that* $\Phi_n$ *is such that* $\tau_{\infty,V}^{(n)}$, $\tau_{V,\infty}^{(n)}$, $\tau_{V,V}^{(n)}$ *and* $\tau_{\infty,\infty}^{(n)}$ *are of order* $o(\frac{n}{\log^3 n})$, $\tau_2^{(n)} := \sup_{\phi \in \Phi_n} \rho_2(\phi) = o(\sqrt{n})$ *and*

$$\int_0^1 \tilde{H}_{\Phi_n}(s)^2 \, \mathrm{d}s = o(\sqrt{n}).$$

*Then*

$$\sup_{\phi \in \Phi_n} |F_n(\phi) - F(\phi)| = \sup_{\phi \in \Phi_n} \left| \frac{1}{\sqrt{n}} E_n(\phi) \right| \xrightarrow{P} 0.$$

# 3. Applications

In this section, we give several examples for the statistic $F_n(\phi)$. In all cases, the results from Section 2 can be applied. As a non-trivial application, the uniform convergence of local Whittle estimates is proved in the next section.

**Example 3.1 (Parametric quasi-likelihood estimation).** In Dahlhaus (2000), it has been shown that

$$\mathcal{L}_n(\theta) := \frac{1}{4\pi} \frac{1}{n} \sum_{t=1}^n \int_{-\pi}^\pi \left\{ \log 4\pi^2 f_\theta\left(\frac{t}{n}, \lambda\right) + \frac{J_n(t/n, \lambda)}{f_\theta(t/n, \lambda)} \right\} \mathrm{d}\lambda \tag{31}$$

is an approximation to $-\log$ Gaussian likelihood of a locally stationary process. The above likelihood is a generalization of the Whittle likelihood (Whittle (1953)) to locally stationary processes. An example for a locally stationary process with finite-dimensional parameter $\theta$ is the tvARMA process from Proposition 2.4 with coefficient functions being polynomials in time. Proving the asymptotic properties of

$$\hat{\theta}_n := \operatorname*{arg\,min}_{\theta \in \Theta} \mathcal{L}_n(\theta)$$

is greatly simplified by using the above properties of the empirical spectral process. We give a brief sketch. Let

$$\mathcal{L}(\theta) := \frac{1}{4\pi} \int_0^1 \int_{-\pi}^\pi \left\{ \log 4\pi^2 f_\theta(u, \lambda) + \frac{f(u, \lambda)}{f_\theta(u, \lambda)} \right\} \mathrm{d}\lambda \, \mathrm{d}u \tag{32}$$

be (up to a constant) the asymptotic Kullback–Leibler divergence between the true process and the fitted model (cf. Dahlhaus (1996), Theorem 3.4 ff) and

$$\theta_0 := \underset{\theta \in \Theta}{\arg\min}\, \mathcal{L}(\theta)$$

the best approximating parameter from $\Theta$ (this is the true parameter if the model is correctly specified). We have

$$\mathcal{L}_n(\theta) - \mathcal{L}(\theta) = \frac{1}{\sqrt{n}} E_n\left(\frac{1}{4\pi} f_\theta^{-1}\right) + R_{\log}(f_\theta)$$

with

$$R_{\log}(f_\theta) := \frac{1}{4\pi} \int_{-\pi}^{\pi} \left[\frac{1}{n} \sum_{t=1}^{n} \log f_\theta\left(\frac{t}{n}, \lambda\right) - \int_0^1 \log f_\theta(u, \lambda)\, du\right] d\lambda. \qquad (33)$$

Thus, ignoring the $R_{\log}$-term, uniform convergence follows from the Glivenko–Cantelli-type Theorem 2.12. The $R_{\log}$-term can be treated as in Dahlhaus and Polonik (2006), Lemma A.2. If $\Theta$ is compact and the minimum $\theta_0$ is unique, this implies consistency of $\widehat\theta_n$. Let $\nabla := (\frac{\partial}{\partial \theta_1}, \ldots, \frac{\partial}{\partial \theta_d})'$. Then

$$\sqrt{n}\nabla\mathcal{L}_n(\theta_0) = E_n\left(\frac{1}{4\pi}\nabla f_\theta^{-1}\right)$$

and

$$\nabla^2 \mathcal{L}_n(\theta) = \frac{1}{\sqrt{n}} E_n\left(\frac{1}{4\pi}\nabla^2 f_\theta^{-1}\right) + \frac{1}{n}\sum_{t=1}^{n}\frac{1}{4\pi}\int_{-\pi}^{\pi}\left(\nabla \log f_\theta\left(\frac{t}{n}, \lambda\right)\right)\left(\nabla \log f_\theta\left(\frac{t}{n}, \lambda\right)\right)' d\lambda.$$

The first term of $\nabla^2 \mathcal{L}_n(\theta)$ converges uniformly to 0, while the second term converges for $\theta \to \theta_0$ to the Fisher information matrix. The usual Taylor expansion then gives a central limit theorem for $\sqrt{n}(\widehat\theta_n - \theta_0)$. For details about the result and examples, we refer to Dahlhaus (2000), Theorem 3.1, the proof of which is greatly simplified by using the above arguments. Furthermore, the present assumptions are weaker. Due to (36) below, the result also covers the misspecified stationary case where the stationary Whittle likelihood is used with a stationary model but the true process is only locally stationary.

Another application of the empirical spectral process is model selection for Whittle estimates. In Van Bellegem and Dahlhaus (2006), a model selection criterion for semi-parametric model selection has been derived. Furthermore, an upper bound for the risk has been proven by using the exponential inequality for the empirical spectral process.

**Example 3.2 (Nonparametric quasi-likelihood estimation).** In Dahlhaus and Polonik (2006), we have considered the corresponding nonparametric estimator

$$\widehat f_n = \underset{g \in \mathcal{F}}{\arg\min}\, \mathcal{L}_n(g)$$

with

$$\mathcal{L}_n(g) = \frac{1}{n} \sum_{t=1}^{n} \frac{1}{4\pi} \int_{-\pi}^{\pi} \left\{ \log g\left(\frac{t}{n}, \lambda\right) + \frac{J_n(t/n, \lambda)}{g(t/n, \lambda)} \right\} \mathrm{d}\lambda, \qquad (34)$$

where the contrast functional is now minimized over an 'infinite-dimensional' target class $\mathcal{F}$ of spectral densities whose complexity is characterized by metric entropy conditions. The optimal rate of convergence has been derived for sieve estimates in the Gaussian case by using a 'peeling device' and 'chaining' together with an exponential inequality similar to the one in Theorem 2.7 (the exponential inequality is stronger due to the additional Gaussian assumption). It is an open problem as to whether the optimal rates of convergence can also be achieved for full nonparametric maximum likelihood estimates or (as in the present paper) without the assumption of Gaussianity.

**Example 3.3 (Stationary processes/model misspecification by stationary models).** We start by showing how several classical results for the stationary case can be obtained from the results above. Let $\phi(u, \lambda) = \tilde{\phi}(\lambda)$ be time invariant. Then

$$F_n(\phi) = \int_{-\pi}^{\pi} \tilde{\phi}(\lambda) \frac{1}{n} \sum_{t=1}^{n} J_n\left(\frac{t}{n}, \lambda\right) \mathrm{d}\lambda. \qquad (35)$$

However, we have

$$
\begin{aligned}
\frac{1}{n} \sum_{t=1}^{n} J_n\left(\frac{t}{n}, \lambda\right) &= \frac{1}{n} \sum_{t=1}^{n} \frac{1}{2\pi} \sum_{1 \le [t+0.5+k/2], [t+0.5-k/2] \le n} X_{[t+0.5+k/2], n} X_{[t+0.5-k/2], n} \exp(-\mathrm{i}\lambda k) \\
&= \frac{1}{2\pi} \sum_{k=-(n-1)}^{n-1} \left( \frac{1}{n} \sum_{t=1}^{n-|k|} X_t X_{t+|k|} \right) \exp(-\mathrm{i}\lambda k) \\
&= \frac{1}{2\pi n} \left| \sum_{s=1}^{n} X_s \exp(-\mathrm{i}\lambda s) \right|^2 = I_n(\lambda),
\end{aligned}
\qquad (36)
$$

where $I_n(\lambda)$ is the classical periodogram. Therefore, $F_n(\phi)$ is the classical spectral mean in the stationary case with the following applications:

(i) $\phi(u, \lambda) = \tilde{\phi}(\lambda) = I_{[0,\mu]}(\lambda)$ gives the empirical spectral measure;

(ii) $\phi(u, \lambda) = \tilde{\phi}(\lambda) = \frac{1}{4\pi} \nabla f_\theta^{-1}(\lambda)$ is the score function of the Whittle likelihood (similar to Example 3.1 above);

(iii) $\phi(u, \lambda) = \tilde{\phi}(\lambda) = \cos \lambda k$ is the empirical covariance estimator of lag $k$.

Theorem 2.5 gives, in all cases, the asymptotic distribution – both in the stationary case and in the misspecified case where the true underlying process is only locally stationary. In case (i) Theorem 2.11 leads to a functional central limit theorem on $C[0, \pi]$ with the supremum-norm. If $\tilde{\phi}(\lambda)$ is a kernel, we obtain a kernel estimate of the spectral density (see the remark below).

***Example 3.4 (Local estimates).*** There is an even larger variety of local estimates
– some of them are listed below. The asymptotic distribution of these estimates is not
covered by Theorem 2.5 since the function $\phi(u, \lambda)$ depends on $n$ in this case. However,
in all cases, the uniform rate of convergence of these estimators may be derived by using
the maximal inequality in Theorem 2.9. A detailed example is given in the next section
where the uniform rate of convergence of local Whittle estimates is derived ((iii) below).

For a short overview, let $k_n(x) = \frac{1}{b_n} K(\frac{x}{b_n})$ be some kernel with bandwidth $b_n$. Then

(i) $\phi(u, \lambda) = k_n(u - u_0) k_n(\lambda - \lambda_0)$  gives an estimator of the time-varying spectral
density $f(u_0, \lambda_0)$;

(ii) $\phi(u, \lambda) = k_n(u - u_0) \cos \lambda k$  gives a local estimator of the covariance function
$c(u_0, k)$;

(iii) $\phi(u, \lambda) = k_n(u - u_0) \frac{1}{4\pi} \nabla f_\theta^{-1}(\lambda)$  is the score function of the local Whittle estima-
tor of the parameter curve $\theta(u_0)$.

# 4. Uniform convergence of local Whittle estimates

We now study kernel estimates for parameter curves of locally stationary processes and
derive uniform consistency from the Glivenko–Cantelli-type Theorem 2.12 (see (39) be-
low) and a uniform rate of convergence from the maximal inequality in Theorem 2.9. We
investigate locally stationary processes where the time varying spectral density is of the
form $f(u, \lambda) = f_{\theta(u)}(\lambda)$ with $\theta(u) \in \Theta \subseteq \mathbf{R}^d$ for all $u \in [0, 1]$. An example is the tvARMA
process from Proposition 2.4.

Let

$$\widehat{\theta}_n(u) := \underset{\theta \in \Theta}{\arg\min}\, \mathcal{L}_n(u, \theta)$$

with

$$\mathcal{L}_n(u, \theta) := \frac{1}{4\pi} \frac{1}{n} \sum_{t=1}^{n} \frac{1}{b_n} K\left(\frac{u - t/n}{b_n}\right) \int_{-\pi}^{\pi} \left\{ \log 4\pi^2 f_\theta(\lambda) + \frac{J_n(t/n, \lambda)}{f_\theta(\lambda)} \right\} \mathrm{d}\lambda. \quad (37)$$

We assume that the kernel $K$ has compact support on $[-\frac{1}{2}, \frac{1}{2}]$ and is of bounded variation
with $\int_{-1/2}^{1/2} x K(x)\, \mathrm{d}x = 0$ and $\int_{-1/2}^{1/2} K(x)\, \mathrm{d}x = 1$. Furthermore, let $b_n \to 0$ and $n b_n \to \infty$
as $n \to \infty$.

In case of a tvAR($p$) process, $\widehat{\theta}_n(u)$ is the solution of the local Yule–Walker equa-
tions: Let $\widehat{c}_n(u, k) := \frac{1}{n} \sum_t \frac{1}{b_n} K(\frac{u - t/n}{b_n}) X_{[t+1/2+k/2], n} X_{[t+1/2-k/2], n}$ (cf. Proposition 5.4),
$C_n(u) = (\widehat{c}_n(u, 1), \ldots, \widehat{c}_n(u, p))'$ and $\Sigma_n(u) = \{\widehat{c}_n(u, i - j)\}_{i, j = 1, \ldots, p}$. If $\widehat{\theta}_n(u) = (\widehat{\alpha}_1(u), \ldots,$
$\widehat{\alpha}_p(u), \widehat{\sigma}^2(u))'$, then it is not difficult to show that

$$(\widehat{\alpha}_1(u), \ldots, \widehat{\alpha}_p(u))' = -\Sigma_n(u)^{-1} C_n(u)$$

and

$$\widehat{\sigma}^2(u) = \widehat{c}_n(u, 0) + \sum_{k=1}^{p} \widehat{\alpha}_k(u)\, \widehat{c}_n(u, k).$$

We now derive a uniform rate of convergence for $\widehat{\theta}_n(u)$. Let $\nabla := (\frac{\partial}{\partial \theta_1}, \ldots, \frac{\partial}{\partial \theta_d})'$, $\|\cdot\|_2$ be the $\ell_2$-norm and $\|A\|_{\mathrm{spec}} := \sup_{x \in \mathbf{C}^n} \frac{\|Ax\|_2}{\|x\|_2}$ be the spectral norm (where $A$ is an $n \times n$ matrix).

**Theorem 4.1.** *Suppose Assumption 2.1 is fulfilled with $\mathbf{E}|\varepsilon_t|^k \leq C_\varepsilon^k$ for all $k \in \mathbf{N}$ and time-varying spectral density $f(u, \lambda) = f_{\theta_0(u)}(\lambda)$. Suppose, further,*

(i) *$\theta$ is identifiable from $f_\theta$ (i.e., $f_\theta(\lambda) = f_{\theta'}(\lambda)$ for all $\lambda$ implies $\theta = \theta'$) and $\theta_0(u)$ lies in the interior of the compact parameter space $\Theta \subseteq \mathbf{R}^d$ for all $u$;*

(ii) *$\theta_0(\cdot)$ is differentiable with Lipschitz continuous derivative;*

(iii) *$f_\theta(\lambda)$ is twice differentiable in $\theta$; $f_\theta^{-1}(\lambda)$ and the components of $\nabla f_\theta(\lambda)$ and $\nabla^2 f_\theta(\lambda)$ are uniformly bounded in $\lambda$ and $\theta$ and uniformly Lipschitz continuous in $\lambda$;*

(iv) *the minimal eigenvalue of $I(\theta) := \frac{1}{4\pi} \int_{-\pi}^{\pi} (\nabla \log f_\theta(\lambda))(\nabla \log f_\theta(\lambda))'\,\mathrm{d}\lambda$ is bounded away from $0$ uniformly in $\theta$.*

*We then have, for $b_n n \gg (\log n)^6$,*

$$\sup_{u \in [b_n/2,\, 1 - b_n/2]} \|\widehat{\theta}_n(u) - \theta_0(u)\|_2 = O_p\left(\frac{1}{\sqrt{b_n n}} + b_n^2\right),$$

*that is, for $b_n \sim n^{-1/5}$, we obtain the uniform rate $O_p(n^{-2/5})$.*

*Remark 4.2.* (i) We conjecture that a similar result also holds in case of model misspecification where the model spectral density $f_{\theta_0(u)}(\lambda)$ is only an approximation to the true spectral density $f(u, \lambda)$.

(ii) It is possible to extend the above result to a wider range of smoothness classes as, for example, in Moulines, Priouret and Roueff (2005). This is not too difficult since only the estimation of the second summand of (41) below needs to be improved in an obvious way (and the kernel $K$ needs to be replaced by a higher order kernel).

(iii) For $\mathrm{tvAR}(p)$ processes, it follows from Theorem 4 in Moulines, Priouret and Roueff (2005) that the above rate is the optimal rate of convergence.

(iv) Of course, the assumption $\mathbf{E}|\varepsilon_t|^k \leq C_\varepsilon^k$ for all $k \in \mathbf{N}$ is restrictive in comparison with what one would expect for standard estimation results based on spectral analysis. It is the price that must be paid for a uniform convergence result with optimal rate and an elegant proof with the maximal exponential inequality.

**Proof of Theorem 4.1.** We begin by noting that the difficult parts of the following proof are handled by using the empirical spectral process and applying Theorems 2.9 and 2.12. We start by proving consistency. We have, with

$$\mathcal{L}(u,\theta) := \frac{1}{4\pi} \int_{-\pi}^{\pi} \left\{ \log 4\pi^2 f_\theta(\lambda) + \frac{f(u,\lambda)}{f_\theta(\lambda)} \right\} d\lambda,$$

$$\mathcal{L}_n(u,\theta) - \mathcal{L}(u,\theta) = \frac{1}{\sqrt{n}} E_n\left( \frac{1}{b_n} K\left( \frac{u-\cdot}{b_n} \right) \otimes \frac{1}{4\pi} f_\theta^{-1} \right)$$

$$+ \frac{1}{4\pi} \int_{-\pi}^{\pi} \int_0^1 \frac{1}{b_n} K\left( \frac{u-v}{b_n} \right) \frac{f(v,\lambda) - f(u,\lambda)}{f_\theta(\lambda)} \, dv \, d\lambda$$

$$+ \frac{1}{4\pi} \int_{-\pi}^{\pi} \log 4\pi^2 f_\theta(\lambda) \, d\lambda \left( \frac{1}{n} \sum_{t=1}^n \frac{1}{b_n} K\left( \frac{u-t/n}{b_n} \right) - 1 \right).$$

We now apply Theorem 2.12 with $\Phi_n = \{ \frac{1}{b_n} K(\frac{u-\cdot}{b_n}) \otimes \frac{1}{4\pi} f_\theta^{-1} \mid u \in [b_n/2, 1 - b_n/2], \theta \in \Theta \}$. It is straightforward to show that $N(\varepsilon, \Phi_n, \rho_2) \leq K / (b_n^{(d+4)/2} \varepsilon^{d+2})$, that is $\int_0^1 \bar{H}_{\Phi_n}(s)^2 \, ds = O((\log b_n)^2)$.

Furthermore, $\tau_{\infty,V}^{(n)}, \tau_{V,\infty}^{(n)}, \tau_{V,V}^{(n)}$ and $\tau_{\infty,\infty}^{(n)}$ are of order $O(b_n^{-1})$ and $\tau_2^{(n)} = \sup_{\phi \in \Phi_n} \rho_2(\phi) = O(b_n^{-1/2})$. Thus, for $b_n n \geq \log^4 n$, Theorem 2.12 implies that

$$\sup_{u \in [b_n/2, 1-b_n/2]} \sup_{\theta \in \Theta} |\mathcal{L}_n(u,\theta) - \mathcal{L}(u,\theta)| \xrightarrow{P} 0. \tag{38}$$

The identifiability condition implies that $\theta_0(u)$ is the unique minimum of $\mathcal{L}(u,\theta)$ for all $u$. By using standard arguments, we can therefore conclude that

$$\sup_{u \in [b_n/2, 1-b_n/2]} \|\widehat{\theta}_n(u) - \theta_0(u)\|_2 \xrightarrow{P} 0. \tag{39}$$

We now derive the rate of convergence by using the maximal inequality. We have, for each $u$,

$$Z_n(u) := \nabla \mathcal{L}_n(u, \widehat{\theta}_n(u)) - \nabla \mathcal{L}_n(u, \theta_0(u)) = \nabla^2 \mathcal{L}_n(u, \bar{\theta}_n(u))(\widehat{\theta}_n(u) - \theta_0(u)) \tag{40}$$

with $|\bar{\theta}_n(u) - \theta_0(u)| \leq |\widehat{\theta}_n(u) - \theta_0(u)|$. The main term on the left-hand side is

$$\frac{\partial}{\partial \theta_j} \mathcal{L}_n(u, \theta_0(u))$$

$$= \frac{1}{\sqrt{n}} E_n(\phi_n) + \int_{-\pi}^{\pi} \int_0^1 \frac{1}{b_n} K\left( \frac{u-v}{b_n} \right) (f_{\theta_0(v)}(\lambda) - f_{\theta_0(u)}(\lambda)) \frac{\partial}{\partial \theta_j} f_\theta^{-1}(\lambda)_{|\theta=\theta_0(u)} \, dv \, d\lambda \tag{41}$$

with $\phi_n(v,\lambda) := \frac{1}{b_n} K(\frac{u-v}{b_n}) \frac{1}{4\pi} \frac{\partial}{\partial \theta_j} f_\theta^{-1}(\lambda)_{|\theta=\theta_0(u)}$. We apply the maximal inequality of Theorem 2.9 to the class $\Phi_n = \{ \frac{1}{b_n} K(\frac{u-\cdot}{b_n}) \otimes \frac{1}{4\pi} \frac{\partial}{\partial \theta_j} f_\theta^{-1} \mid u \in [b_n/2, 1-b_n/2], \theta \in \Theta \}$. Again, we can show $N(\varepsilon, \Phi_n, \rho_2) = K / (b_n^{(d+4)/2} \varepsilon^{d+2})$, that is, $\int_0^1 \tilde{H}_{\Phi_n}(s)^2 \, ds = O((\log b_n)^2)$. Furthermore, $\tau_{\infty,V}^{(n)}$, $\tau_{V,\infty}^{(n)}$, $\tau_{V,V}^{(n)}$ and $\tau_{\infty,\infty}^{(n)}$ are of order $O(b_n^{-1})$ and $\tau_2^{(n)} = \sup_{\phi \in \Phi_n} \rho_2(\phi) \sim b_n^{-1/2}$. We now apply Theorem 2.9 with $\eta = \tau_2 \delta$ for arbitrary $\delta$. If $b_n n \gg \log^6 n$, then the conditions (26) and (27) are fulfilled and we obtain

$$P\left( \sup_{\phi \in \Phi_n} |E_n(\phi)| > \tau_2 \delta, \; B_n \right) \leq 3c_1 \exp\left\{ -\frac{c_2}{4} \sqrt{\delta} \right\}$$

and, as a consequence,

$$\sup_{\phi \in \Phi_n} \left\| \frac{1}{\sqrt{n}} E_n(\phi) \right\| = O_p\left( \frac{1}{\sqrt{b_n n}} \right).$$

The smoothness conditions (ii) and (iii) imply that the second summand of (41) can be uniformly bounded by $O(b_n^2)$, that is, we obtain

$$\sup_{u \in [b_n/2, 1-b_n/2]} \|\nabla \mathcal{L}_n(u, \theta_0(u))\|_2 = O_p\left( \frac{1}{\sqrt{b_n n}} + b_n^2 \right).$$

If $\hat{\theta}_n(u)$ lies on the boundary of $\Theta$ for some $u$, then $\|\hat{\theta}_n(u) - \theta(u)\|_2 \geq \kappa$ for some $\kappa > 0$ and, by (39),

$$P\left( \sup_{u \in [b_n/2,\, 1-b_n/2]} \|\nabla \mathcal{L}_n(u, \hat{\theta}_n(u))\|_2 > \delta \frac{1}{\sqrt{b_n n}} \right)$$

$$\leq P\left( \sup_{u \in [b_n/2, 1-b_n/2]} \|\nabla \mathcal{L}_n(u, \hat{\theta}_n(u))\|_2 > 0 \right)$$

$$\leq P\left( \sup_{u \in [b_n/2, 1-b_n/2]} \|\hat{\theta}_n(u) - \theta_0(u)\|_2 \geq \kappa \right) \to 0,$$

implying that $\sup_{u \in [b_n/2, 1-b_n/2]} \|Z_n(u)\|_2 = O_p(\frac{1}{\sqrt{b_n n}})$. In order to obtain the assertion of the theorem from (40), we now prove that the minimal eigenvalue of the matrix $\nabla^2 \mathcal{L}_n(u, \bar{\theta}_n(u))$ is bounded away from 0 uniformly in $\theta$ in probability. We have

$$\nabla^2 \mathcal{L}_n(u, \theta) = \frac{1}{\sqrt{n}} E_n\left( \frac{1}{b_n} K\left( \frac{u - \cdot}{b_n} \right) \otimes \frac{1}{4\pi} \nabla^2 f_\theta^{-1} \right)$$

$$+ \frac{1}{n} \sum_{t=1}^n \frac{1}{b_n} K\left( \frac{u - t/n}{b_n} \right) \frac{1}{4\pi} \int_{-\pi}^{\pi} (\nabla \log f_\theta(\lambda))(\nabla \log f_\theta(\lambda))' \, d\lambda.$$

Since $f_\theta$ is twice differentiable in $\theta$ with Lipschitz continuous second derivative in $\lambda$ we obtain, exactly as above from Theorem 2.12 for $b_n n \geq \log^4 n$ and $i, j = 1, \ldots, d$,

$$\sup_{u \in [b_n/2,\, 1-b_n/2]} \sup_{\theta \in \Theta} \left| \frac{1}{\sqrt{n}} E_n\left( \frac{1}{b_n} K\left( \frac{u - \cdot}{b_n} \right) \otimes \frac{1}{4\pi} \frac{\partial^2}{\partial \theta_i \partial \theta_j} f_\theta^{-1} \right) \right| \xrightarrow{P} 0.$$

Therefore, also,

$$\sup_{u \in [b_n/2, 1-b_n/2]} \left\| \frac{1}{\sqrt{n}} E_n \left( \frac{1}{b_n} K \left( \frac{u - \cdot}{b_n} \right) \otimes \frac{1}{4\pi} \left( \nabla^2 f_\theta^{-1} \right)_{|\theta = \bar{\theta}_n(u)} \right) \right\|_{\text{spec}} \xrightarrow{P} 0.$$

Since the minimal eigenvalue of $I(\theta) := \frac{1}{4\pi} \int_{-\pi}^{\pi} (\nabla \log f_\theta(\lambda))(\nabla \log f_\theta(\lambda))' \, d\lambda$ is bounded from below by $\lambda_{\min}(I) > 0$ uniformly in $\theta$, this implies that

$$\mathbf{P} \left( \sup_{u \in [b_n/2, \, 1-b_n/2]} \|\nabla^2 \mathcal{L}_n(u, \bar{\theta}_n(u))^{-1}\|_{\text{spec}} \le \frac{2}{\lambda_{\min}(I)} \right) \to 1.$$

Since

$$\|\widehat{\theta}_n(u) - \theta_0(u)\|_2 \le \|\nabla^2 \mathcal{L}_n(u, \bar{\theta}_n(u))^{-1}\|_{\text{spec}} \|Z_n(u)\|_2,$$

this implies the result.                                                                      $\square$

## 5. Proofs: CLT and exponential inequality

In this section, we provide the proofs for the results of Section 2. In particular, we derive the asymptotic behavior of the moments of the empirical spectral process.

First, we extend the definitions of Section 2 to tapered data $X_{t,n}^{(h_n)} = h_n(\frac{t}{n}) \cdot X_{t,n}$, where $h_n : (0,1] \to [0, \infty)$ is a data taper (with $h_n(\cdot) = I_{(0,1]}(\cdot)$ being the non-tapered case of Section 2). This is done for three reasons:

(i) The main reason is that all proofs are greatly simplified since the data taper now automatically takes care of the range of summation ($h_n(t/n)$ is zero for all $t$ outside the observation domain $\{1, \ldots, n\}$). The consideration of arbitrary tapers $h_n$ instead of the 'no-taper' $I_{(0,1]}$ does not introduce any extra technical complexity at all.

(ii) By using tapers which are different from 0 only on a segment of the observation domain, one may construct localized estimators. This is of great importance, although we will not discuss it in the present paper.

(iii) The use of a data taper in the periodogram is standard for stationary time series. It leads to a better small-sample performance in the presence of strong peaks in the spectrum. It may turn out that this also holds in the present situation (which requires further investigation).

As before, the empirical spectral process is defined by $E_n^{(h_n)}(\phi) := \sqrt{n} \, (F_n^{(h_n)}(\phi) - F^{(h_n)}(\phi))$, where

$$F^{(h_n)}(\phi) := \int_0^1 h_n^2(u) \int_{-\pi}^{\pi} \phi(u, \lambda) f(u, \lambda) \, d\lambda \, du \tag{42}$$

and

$$F_n^{(h_n)}(\phi) := \frac{1}{n} \sum_{t=1}^n \int_{-\pi}^{\pi} \phi \left( \frac{t}{n}, \lambda \right) J_n^{(h_n)} \left( \frac{t}{n}, \lambda \right) d\lambda, \tag{43}$$

now with the tapered pre-periodogram

$$J_n^{(h_n)}\left(\frac{t}{n}, \lambda\right) = \frac{1}{2\pi} \sum_{k:1\leq[t+1/2\pm k/2]\leq n} X_{[t+1/2+k/2],n}^{(h_n)} X_{[t+1/2-k/2],n}^{(h_n)} \exp(-\mathrm{i}\lambda k). \qquad (44)$$

We mention that in some cases, a rescaling may be necessary for $J_n^{(h_n)}(u, \lambda)$ to become a (pre-) estimate of $f(u, \lambda)$.

**Throughout this appendix the superscript $(h_n)$ will be dropped in many situations for notational convenience, that is, we will use $F_n(\phi)$, $F(\phi)$ and $E_n(\phi)$.**

**Assumption 5.1.** *The data taper $h_n : (0, 1] \to [0, \infty)$ fulfills $\sup_n V(h_n) \leq C$ and $\sup_{u,n} h_n(u) \leq C$ for some $C < \infty$. Furthermore, $\log h_n(\cdot)$ is concave.*

The assumption that $\log h_n(\cdot)$ is concave is very mild (note that even $\log(x^m)$ is concave). We need the following notation. With

$$\hat{\phi}(u, j) := \int_{-\pi}^{\pi} \phi(u, \lambda) \exp(\mathrm{i}\lambda j) \, \mathrm{d}\lambda, \qquad (45)$$

we define

$$\rho_\infty(\phi) := \sum_{j=-\infty}^{\infty} \sup_u |\hat{\phi}(u, j)| \quad \text{and} \quad v_\Sigma(\phi) := \sum_{j=-\infty}^{\infty} V(\hat{\phi}(\cdot, j)). \qquad (46)$$

We mention that

$$\sup_{u,\lambda} |\phi(u, \lambda)| \leq \frac{1}{2\pi} \rho_\infty(\phi) \quad \text{and} \quad \rho_2(\phi) \leq \frac{1}{\sqrt{2\pi}} \rho_\infty(\phi).$$

The idea now is to prove the CLT in Theorem 2.5 by the convergence of all cumulants. The convergence of the cumulants is derived below under the assumptions $\rho_\infty(\phi) < \infty$ and $v_\Sigma(\phi) < \infty$. Unfortunately, these assumptions are not fulfilled for functions of bounded variation as assumed in Theorem 2.5. Therefore, its proof also uses certain approximation arguments (cf. proof of Theorem 2.5). The following CLT is a by-product which follows immediately from the cumulant calculations below. It is of independent interest since the result does not follow from Theorem 2.5 (the condition $\rho_\infty(\phi) < \infty$ does not imply bounded variation in the $\lambda$-direction; furthermore, the conditions may be easier to check in some situations).

**Assumption 5.2.** *Suppose $\phi : [0, 1] \times [-\pi, \pi] \to \mathbf{R}$ fulfills $\rho_\infty(\phi) < \infty$ and $v_\Sigma(\phi) < \infty$.*

**Theorem 5.3.** *Suppose Assumptions 2.1, 5.1 and 5.2 hold with a data taper $h$ independent of $n$. Then*

$$(E_n^{(h)}(\phi_j))_{j=1,\ldots,k} \overset{\mathcal{D}}{\to} (E^{(h)}(\phi_j))_{j=1,\ldots,k},$$

where $(E^{(h)}(\phi_j))_{j=1,\dots,k}$ is a Gaussian random vector with mean $0$ and $\operatorname{cov}(E^{(h)}(\phi_j),$ $E^{(h)}(\phi_k)) = c_E^{(h)}(\phi_j, \phi_k)$ with

$$c_E^{(h)}(\phi_j, \phi_k) := 2\pi \int_0^1 h^4(u) \int_{-\pi}^{\pi} \phi_j(u, \lambda) \left[\phi_k(u, \lambda) + \phi_k(u, -\lambda)\right] f^2(u, \lambda) \, \mathrm{d}\lambda \, \mathrm{d}u$$

$$+ \kappa_4 \int_0^1 h^4(u) \left(\int_{-\pi}^{\pi} \phi_j(u, \lambda_1) f(u, \lambda_1) \, \mathrm{d}\lambda_1\right) \left(\int_{-\pi}^{\pi} \phi_k(u, \lambda_2) f(u, \lambda_2) \, \mathrm{d}\lambda_2\right) \mathrm{d}u.$$

**Proof.** The result follows from the convergence of all cumulants which is proved in Lemma 5.5(ii), Lemma 5.6(ii) and Lemma 5.7(iii). $\qquad\square$

In the same way, one may arrive at a central limit theorem with $(\|h_n\|_2^{-1} E_n(\phi_j))_{j=1,\dots,k}$ for $h_n$ dependent on $n$ (under additional assumptions), for example, for $h_n(\cdot) = I_{[u_0 - b_n/2,\, u_0 - b_n/2]}(\cdot)$ (segment estimate). This will be studied in future work.

For the following proofs, we first need a result on the behavior (decay) of the covariances of the process. The case $h_n(\cdot) = I_{(0,1]}(\cdot)$ gives the results for the ordinary covariances.

**Proposition 5.4.** *Suppose Assumptions 2.1 and 5.1 hold. We then have, for all $k, k_1, k_2 \in \mathbf{Z}$ with some $K$ independent of $k, k_1, k_2$ and $n$,*

$$\sup_t |\operatorname{cov}(X_{t,n}^{(h_n)}, X_{t+k,n}^{(h_n)})| \le \frac{K}{\ell(k)}, \tag{47}$$

$$\sup_u |c(u, k)| \le \frac{K}{\ell(k)}, \tag{48}$$

$$\sum_{t=1}^n \left| \operatorname{cov}(X_{t+k_1,n}^{(h_n)} X_{t-k_2,n}^{(h_n)}) - h_n\left(\frac{t}{n}\right)^2 c\left(\frac{t}{n}, k_1 + k_2\right) \right| \le K \left(1 + \frac{\min\{|k_1|, n\}}{\ell(k_1 + k_2)}\right), \tag{49}$$

$$V(c(\cdot, k)) \le \frac{K}{\ell(k)}. \tag{50}$$

**Proof.** From (11) and (4), we have

$$c(u, k) = \sum_{j=-\infty}^{\infty} a(u, j+k) \, a(u, j) \quad \text{with } \sup_u |a(u, j)| \le \frac{K}{\ell(j)}.$$

Therefore, (48) follows from the relation

$$\sum_{j=-\infty}^{\infty} \frac{1}{\ell(k+j)} \frac{1}{\ell(j)} \le \frac{K}{\ell(k)} \tag{51}$$

which is easily established. Furthermore, we have, with $k = k_1 + k_2$,

$$\text{cov}(X_{t+k_1,n}^{(h_n)}, X_{t-k_2,n}^{(h_n)}) = h_n\left(\frac{t+k_1}{n}\right) h_n\left(\frac{t-k_2}{n}\right) \sum_{j=-\infty}^{\infty} a_{t+k_1,n}(j+k)\, a_{t-k_2,n}(j).$$

For $k_1 = k$ and $k_2 = 0$, this gives (47) by using (3). Replacing $h_n(\frac{t+k_1}{n})$, $h_n(\frac{t-k_2}{n})$, $a_{t+k_1,n}(j+k)$ and $a_{t-k_2,n}(j)$ by $h_n(\frac{t}{n})$, $h_n(\frac{t}{n})$, $a(\frac{t}{n}, j+k)$ and $a(\frac{t}{n}, j)$, respectively, gives (49). For example, the last replacement step has, with (5) and (6), the upper bound

$$\sum_{t=1}^{n} \left| h_n\left(\frac{t}{n}\right) h_n\left(\frac{t}{n}\right) \sum_{j=-\infty}^{\infty} a\left(\frac{t}{n}, j+k\right)\left(a_{t-k_2,n}(j) - a\left(\frac{t}{n}, j\right)\right)\right|$$

$$\leq K \sum_{j=-\infty}^{\infty} \frac{1}{\ell(j+k)} \sum_{t=1}^{n} \left| a_{t-k_2,n}(j) - a\left(\frac{t}{n}, j\right)\right|$$

$$\leq K \sum_{j=-\infty}^{\infty} \frac{1}{\ell(j+k)}\left(1 + \frac{|k_2|}{\ell(j)}\right)$$

$$= K\left(1 + \frac{|k_2|}{\ell(k)}\right).$$

Since $|k_2| \leq |k| + |k_1|$, we obtain (49) for $|k_1| \leq n$. For $|k_1| > n$, the result follows since then $\text{cov}(X_{t+k_1,n}^{(h_n)}, X_{t-k_2,n}^{(h_n)})$ is equal to 0. (50) is obtained in the same way. $\qquad\square$

A trick which greatly simplifies the following proofs is to set $a_{t,n}(j) = 0$ for $t \notin \{1, \ldots, n\}$ and $j \in \mathbf{Z}$, $a(u, j) = 0$ for $u \notin (0, 1]$ and $j \in \mathbf{Z}$, $\phi(u, \lambda) = 0$ for $u \notin (0, 1]$ and $\lambda \in [-\pi, \pi]$ and $h_n(u) = 0$ for $u \notin (0, 1]$. With this convention, (3)–(6), (47)–(50) continue to hold for $u \in \mathbf{R}$, $t \in \mathbf{Z}$, $V(f)$ now denoting the total variation over $\mathbf{R}$ and $t$ in (5) and (49) ranging from $-\infty$ to $\infty$. Furthermore, the summation range of $k$ in (16) and $t$ in (15) can be extended from $-\infty$ to $\infty$. Therefore, all summation ranges are from $-\infty$ to $\infty$ in the following proofs unless otherwise indicated.

We also set $\tilde{a}(j) = \sup_u |a(u, j)|$, $\tilde{\phi}(j) = \max\{\sup_u |\hat{\phi}(u, j)|, \sup_u |\hat{\phi}(u, -j)|\}$ and $\tilde{c}(j) = \sup_u |c(u, j)|$.

**Lemma 5.5.** (i) *Suppose Assumptions 2.1 and 5.1 hold and $\phi : [0, 1] \times [-\pi, \pi] \to \mathbf{R}$ is a function possibly depending on $n$. We then have, with $K > 0$,*

$$|\mathbf{E}F_n^{(h_n)}(\phi) - F^{(h_n)}(\phi)| \leq \frac{K}{n} \sum_{|k| \leq n} \tilde{\phi}(k) + K \sum_{|k| > n} \tilde{\phi}(k)\frac{1}{\ell(k)} + \frac{K}{n} \sum_k V\left(c(\cdot, k)\right)\frac{1}{\ell(k)}.$$

(ii) *If, in addition, $\phi$ is independent of $n$ and fulfills Assumption 5.2, then*

$$\mathbf{E}E_n^{(h_n)}(\phi) = O(n^{-1/2}).$$

**Proof.** (i) We have

$$F_n^{(h_n)}(\phi) = \frac{1}{2\pi n} \sum_{t=1}^{n} \sum_k \hat{\phi}\left(\frac{t}{n}, -k\right) X_{[t+1/2+k/2],n}^{(h_n)} X_{[t+1/2-k/2],n}^{(h_n)}. \tag{52}$$

We therefore obtain from Proposition 5.4

$$\mathbf{E} F_n^{(h_n)}(\phi) = \frac{1}{2\pi n} \sum_{t,|k|\leq n} \hat{\phi}\left(\frac{t}{n}, -k\right) \mathrm{cov}(X_{[t+1/2+k/2],n}^{(h_n)}, X_{[t+1/2-k/2],n}^{(h_n)})$$

$$= \frac{1}{2\pi n} \sum_{t,k} h_n^2\left(\frac{t}{n}\right) \hat{\phi}\left(\frac{t}{n}, -k\right) c\left(\frac{t}{n}, k\right) + R \tag{53}$$

with

$$|R| \leq \frac{K}{n} \sum_{|k|\leq n} \tilde{\phi}(k) \left[1 + \frac{\min(|k|, n)}{\ell(k)}\right] + K \sum_{|k|>n} \tilde{\phi}(k) \frac{1}{\ell(k)} \leq \frac{K}{n} \rho_\infty(\phi). \tag{54}$$

Furthermore,

$$\left| \frac{1}{2\pi n} \sum_{t,k} h_n^2\left(\frac{t}{n}\right) \hat{\phi}\left(\frac{t}{n}, -k\right) c\left(\frac{t}{n}, k\right) - F(\phi) \right|$$

$$\leq \left| \frac{1}{2\pi} \sum_{t,k} \int_0^{1/n} \left[ h_n^2\left(\frac{t}{n}\right) \hat{\phi}\left(\frac{t}{n}, -k\right) c\left(\frac{t}{n}, k\right) \right.\right.$$

$$\left.\left. - h_n^2\left(\frac{t-1}{n} + x\right) \hat{\phi}\left(\frac{t-1}{n} + x, -k\right) c\left(\frac{t-1}{n} + x, k\right) \right] \mathrm{d}x \right| \tag{55}$$

$$\leq \frac{K}{2\pi n} \sum_k \left[ V(\hat{\phi}(\cdot, -k)) \tilde{c}(k) + \tilde{\phi}(k) V(c(\cdot, k)) + \tilde{\phi}(k) \tilde{c}(k) \right]$$

leading to the result. (ii) follows immediately. $\qquad\square$

**Lemma 5.6.** (i) *Suppose Assumptions 2.1 and 5.1 hold and $\phi_1, \phi_2 : [0,1] \times [-\pi, \pi] \to \mathbf{R}$ are functions possibly depending on $n$. We then have*

$$\mathrm{cov}(E_n^{(h_n)}(\phi_1), E_n^{(h_n)}(\phi_2)) = c_E^{(h_n)}(\phi_1, \phi_2) + R_n$$

*with*

$$|R_n| \leq \frac{K}{n} \sum_{k_1, k_2} \tilde{\phi}_1(k_1) \tilde{\phi}_2(k_2) \left[ 1 + \frac{\min\{|k_1|, n\}}{\ell(k_1 + k_2)} \right]$$

$$+ \frac{K}{n} \sum_{k_1, k_2, k_3} [\tilde{\phi}_1(k_1) V(\hat{\phi}_2(\cdot, k_2)) V(\hat{\phi}_1(\cdot, k_1)) \tilde{\phi}_2(k_2)] \frac{\min\{|k_1| + |k_2| + |k_3|, n\}}{\ell(k_3)\ell(k_1 + k_2 + k_3)}$$

$$+ \frac{K}{n} \sum_{k_1,k_2} [\tilde{\phi}_1(k_1) V(\hat{\phi}_2(\cdot,k_2)) + V(\hat{\phi}_1(\cdot,k_1)) \tilde{\phi}_2(k_2)] \left[ \frac{1}{\ell(k_1)} + \frac{1}{\ell(k_2)} \right],$$

*where the last term can be omitted if $X_{t,n}$ is Gaussian.*

(ii) *If $\phi_1$ and $\phi_2$ are independent of $n$ and fulfill Assumption 5.2, then $R_n = o(1)$.*

**Proof.** (i) We have, with (52),

$$\begin{aligned}
\mathrm{cov}&(E_n^{(h_n)}(\phi_1), E_n^{(h_n)}(\phi_2)) \\
&= n\, \mathrm{cov}(F_n^{(h_n)}(\phi_1), F_n^{(h_n)}(\phi_2)) \\
&= \frac{1}{(2\pi)^2 n} \sum_{t_1,t_2,k_1,k_2} \hat{\phi}_1\left(\frac{t_1}{n}, -k_1\right) \hat{\phi}_2\left(\frac{t_2}{n}, -k_2\right) \\
&\quad \times [\mathrm{cov}(X_{[t_1+1/2+k_1/2],n}^{(h_n)}, X_{[t_2+1/2+k_2/2],n}^{(h_n)}) \\
&\qquad \times \mathrm{cov}(X_{[t_1+1/2-k_1/2],n}^{(h_n)}, X_{[t_2+1/2-k_2/2],n}^{(h_n)}) \\
&\qquad + \mathrm{cov}(X_{[t_1+1/2+k_1/2],n}^{(h_n)}, X_{[t_2+1/2-k_2/2],n}^{(h_n)}) \\
&\qquad \times \mathrm{cov}(X_{[t_1+1/2-k_1/2],n}^{(h_n)}, X_{[t_2+1/2+k_2/2],n}^{(h_n)}) \\
&\qquad + \mathrm{cum}(X_{[t_1+1/2+k_1/2],n}^{(h_n)}, X_{[t_1+1/2-k_1/2],n}^{(h_n)}, \\
&\qquad\qquad X_{[t_2+1/2+k_2/2],n}^{(h_n)}, X_{[t_2+1/2-k_2/2],n}^{(h_n)})].
\end{aligned} \tag{56}$$

Let $k_3 := t_1 - t_2 + [k_1/2 + 1/2] - [k_2/2 + 1/2]$. By using Proposition 5.4, we replace the first summand in $[\ldots]$ by $h_n(\frac{t_1}{n})^4 c(\frac{t_1}{n}, k_3) c(\frac{t_1}{n}, k_3 + k_2 - k_1)$. The remainder can be bounded by

$$\begin{aligned}
\frac{K}{n} \sum_{k_1,k_2,k_3} \Big[ &\tilde{\phi}_1(k_1)\tilde{\phi}_2(k_2) \Big\{ 1 + \frac{\min\{|k_1|, n\}}{\ell(k_3)} \Big\} \frac{1}{\ell(k_3+k_2-k_1)} \\
&+ \tilde{\phi}_1(k_1)\tilde{\phi}_2(k_2) \frac{1}{\ell(k_3)} \Big\{ 1 + \frac{\min\{|k_1|,n\}}{\ell(k_3+k_2-k_1)} \Big\} \Big].
\end{aligned}$$

(51) implies that this is bounded as asserted. Therefore, the first term is equal to

$$\begin{aligned}
\frac{1}{(2\pi)^2 n} \sum_{t_1,k_1,k_2,k_3} h_n&\left(\frac{t_1}{n}\right)^4 \hat{\phi}_1\left(\frac{t_1}{n}, -k_1\right) \hat{\phi}_2\left(\frac{t_1+k_o}{n}, -k_2\right) \\
&\times c\left(\frac{t_1}{n}, k_3\right) c\left(\frac{t_1}{n}, k_3+k_2-k_1\right) + R_n,
\end{aligned} \tag{57}$$

where $k_o = -k_3 + [k_1/2 + 1/2] - [k_2/2 + 1/2]$. Replacing $\hat{\phi}_2(\frac{t+k_o}{n}, -k_2)$ by $\hat{\phi}_2(\frac{t_1}{n}, -k_2)$ yields the error term

$$\frac{K}{n} \sum_{k_1, k_2, k_3} \tilde{\phi}(k_1)\tilde{c}(k_3)\tilde{c}(k_3 + k_2 - k_1) \sum_t \left| \hat{\phi}_2\left(\frac{t+k_o}{n}, -k_2\right) - \hat{\phi}_2\left(\frac{t}{n}, -k_2\right) \right|$$

which is also bounded as claimed. As in (55), we now replace the $\frac{1}{n}\sum_{t_1}$ sum in (58) (with $k_o = 0$) by the integral over $[0, 1]$ with the same replacement error. Direct calculation (or repeated application of Parseval's equality) yields

$$\frac{1}{(2\pi)^2} \int_0^1 h_n^4(u) \sum_{k_1, k_2, k_3} \hat{\phi}_1(u, -k_1)\, \hat{\phi}_2(u, -k_2)\, c(u, k_3)\, c(u, k_3 + k_2 - k_1)\, \mathrm{d}u$$

$$= 2\pi \int_0^1 h_n^4(u) \int_{-\pi}^{\pi} \phi_1(u, \lambda)\phi_2(u, -\lambda)f(u, \lambda)^2 \, \mathrm{d}\lambda \, \mathrm{d}u.$$

The second term in (56) is treated in the same way. With the representation $X_{t,n} = \sum_{j=-\infty}^{\infty} a_{t,n}(t-j)\varepsilon_j$ and the abbreviations $t_\nu^+ = t_\nu^+(t_\nu, k_\nu) = [t_\nu + 1/2 + k_\nu/2]$, $t_\nu^- = t_\nu^-(t_\nu, k_\nu) = [t_\nu + 1/2 - k_\nu/2]$, the third term is equal to

$$\frac{\kappa_4}{(2\pi)^2 n} \sum_{t_1, t_2, k_1, k_2} \hat{\phi}_1\left(\frac{t_1}{n}, -k_1\right) \hat{\phi}_2\left(\frac{t_2}{n}, -k_2\right)$$

$$\times \sum_i h_n\left(\frac{t_1^+}{n}\right) h_n\left(\frac{t_1^-}{n}\right) h_n\left(\frac{t_2^+}{n}\right) h_n\left(\frac{t_2^-}{n}\right)$$

$$\times a_{t_1^+, n}(t_1^+ - i)\, a_{t_1^-, n}(t_1^- - i)\, a_{t_2^+, n}(t_2^+ - i)\, a_{t_2^-, n}(t_2^- - i).$$

By using Assumption 2.1 and (51), we now replace this by

$$\frac{\kappa_4}{(2\pi)^2 n} \sum_{t_1, t_2, k_1, k_2} \hat{\phi}_1\left(\frac{t_1}{n}, -k_1\right) \hat{\phi}_2\left(\frac{t_2}{n}, -k_2\right)$$

$$\times \sum_i h_n\left(\frac{t_1}{n}\right)^2 h_n\left(\frac{t_2}{n}\right)^2 a\left(\frac{t_1}{n}, t_1^+ - i\right) \tag{58}$$

$$\times a\left(\frac{t_1}{n}, t_1^- - i\right) a\left(\frac{t_2}{n}, t_2^+ - i\right) a\left(\frac{t_2}{n}, t_2^- - i\right)$$

with replacement error $K n^{-1} \sum_{k_1, k_2} \tilde{\phi}_1(k_1)\tilde{\phi}_2(k_2)$. We now replace the term $a(\frac{t_2}{n}, t_2^- - i)$ in the above expression by $a(\frac{t_1}{n}, t_2^- - i)$, leading, with the substitutions $d = t_2 - t_1$ and $j = i - t_1$, to a replacement error of

$$\frac{K}{n} \sum_{k_1, k_2} \tilde{\phi}_1(k_1)\tilde{\phi}_2(k_2)$$

$$\times \sum_{d,j} \frac{1}{\ell([\frac{1}{2}+\frac{k_1}{2}]-j)} \frac{1}{\ell([\frac{1}{2}-\frac{k_1}{2}]-j)} \frac{1}{\ell([d+\frac{1}{2}+\frac{k_2}{2}]-j)}$$

$$\times \sum_{t_1} \left| a\left(\frac{t_1+d}{n}, \left[d+\frac{1}{2}-\frac{k_2}{2}\right]-j\right) - a\left(\frac{t_1}{n}, \left[d+\frac{1}{2}-\frac{k_2}{2}\right]-j\right)\right|.$$

The last sum is bounded by $|d| V(a(\cdot, [d+\frac{1}{2}-\frac{k_2}{2}]-j))$, leading to the upper bound

$$\frac{K}{n} \sum_{k_1, k_2} \tilde{\phi}_1(k_1) \tilde{\phi}_2(k_2)$$

$$\times \sum_{d,j} \frac{|[d+\frac{1}{2}+\frac{k_2}{2}]-j| + |[d+\frac{1}{2}-\frac{k_2}{2}]-j| + |[\frac{1}{2}+\frac{k_1}{2}]-j| + |[\frac{1}{2}-\frac{k_1}{2}]-j|}{\ell([\frac{1}{2}+\frac{k_1}{2}]-j)\,\ell([\frac{1}{2}-\frac{k_1}{2}]-j)\,\ell([d+\frac{1}{2}+\frac{k_2}{2}]-j)\,\ell([d+\frac{1}{2}-\frac{k_2}{2}]-j)}$$

$$\leq \frac{K}{n} \sum_{k_1, k_2} \tilde{\phi}_1(k_1) \tilde{\phi}_2(k_2)$$

for the replacement error. In the same way, we replace (59) by

$$\frac{\kappa_4}{(2\pi)^2 n} \sum_{t_1, t_2, k_1, k_2} \hat{\phi}_1\left(\frac{t_1}{n}, -k_1\right) \hat{\phi}_2\left(\frac{t_1}{n}, -k_2\right)$$

$$\times \sum_i h_n\left(\frac{t_1}{n}\right)^4 a\left(\frac{t_1}{n}, t_1^+ - i\right)$$

$$\times a\left(\frac{t_1}{n}, t_1^- - i\right) a\left(\frac{t_1}{n}, t_2^+ - i\right) a\left(\frac{t_1}{n}, t_2^- - i\right) \tag{59}$$

$$= \frac{\kappa_4}{(2\pi)^2 n} \sum_{t_1} h_n\left(\frac{t_1}{n}\right)^4 \sum_{k_1, k_2} \hat{\phi}_1\left(\frac{t_1}{n}, -k_1\right) \hat{\phi}_2\left(\frac{t_1}{n}, -k_2\right) c\left(\frac{t_1}{n}, k_1\right) c\left(\frac{t_1}{n}, k_2\right)$$

with replacement error $K n^{-1} \sum_{k_1, k_2} \tilde{\phi}_1(k_1)[\tilde{\phi}_2(k_2) + V(\hat{\phi}_2(\cdot, k_2))][\frac{1}{\ell(k_1)} + \frac{1}{\ell(k_2)}]$. As in (55), we now replace the $\frac{1}{n}\sum_{t_1}$ sum by the integral over $[0, 1]$. Application of Parseval's equality gives the final form of the fourth order cumulant term.

(ii) Considering the cases $|k| \leq \sqrt{n}$ and $|k| > \sqrt{n}$ separately shows that

$$\frac{1}{n} \sum_k \min\{|k|, n\} \tilde{\phi}_i(k) = o(1) \quad \text{and} \quad \frac{1}{n} \sum_k \min\{|k|, n\} \frac{1}{\ell(k)} = o(1). \tag{60}$$

This implies that the first term of $R_n$ tends to zero. Since

$$\min\{|k_1| + |k_2| + |k_3|, n\}$$

$$\leq 2\min\{|k_1 + k_2 + k_3|, n\} + 2\min\{|k_1|, n\} + 2\min\{|k_3|, n\}$$

and $|V(\hat{\phi}_i(\cdot, k))| \leq K$, the third term of $R_n$ also tends to zero. $\qquad \square$

We now set

$$\rho_{2,n}^{(h_n)}(\phi) := \left( \frac{1}{n} \sum_{t=1}^{n} h_n \left( \frac{[t+1/2]}{n} \right)^4 \int_{-\pi}^{\pi} \phi \left( \frac{t}{n}, \lambda \right)^2 d\lambda \right)^{1/2}. \tag{61}$$

Note that $\rho_{2,n}^{(h_n)}(\phi) = \rho_{2,n}(\phi)$ in the non-tapered case where $h_n(\cdot) = I_{(0,1]}(\cdot)$.

**Lemma 5.7.** *Suppose Assumptions 2.1 and 5.1 hold and $\phi_1, \ldots, \phi_\ell : [0,1] \times [-\pi, \pi] \to \mathbf{R}$ are functions possibly depending on $n$.*

(i) *If $\ell \geq 2$, then*

$$|\mathrm{cum}(E_n^{(h_n)}(\phi_1), \ldots, E_n^{(h_n)}(\phi_\ell))| \leq K n^{1-\ell/2} \rho_{2,n}^{(h_n)}(\phi_1) \, \rho_{2,n}^{(h_n)}(\phi_2) \prod_{j=3}^{\ell} \rho_\infty(\phi_j)$$

*with a constant $K$ independent of $n$.*

(ii) *If $\ell \geq 2$ and, in addition, $\mathbf{E}|\varepsilon_t|^k \leq C_\varepsilon^k$ for all $k \in \mathbf{N}$ for the $\varepsilon_t$ from Assumption 2.1, then*

$$|\mathrm{cum}(E_n^{(h_n)}(\phi_1), \ldots, E_n^{(h_n)}(\phi_\ell))| \leq K^\ell (2\ell)! \prod_{j=1}^{\ell} \rho_{2,n}^{(h_n)}(\phi_j) \tag{62}$$

*with a constant $K$ independent of $n$ and $\ell$.*

(iii) *If $\ell \geq 3$ and $\phi_1, \ldots, \phi_\ell$ are independent of $n$ and fulfill Assumption 5.2, then*

$$|\mathrm{cum}(E_n^{(h_n)}(\phi_1), \ldots, E_n^{(h_n)}(\phi_\ell))| = O(n^{1-\ell/2}).$$

**Proof.** (i) We have, with (52),

$$\mathrm{cum}(E_n^{(h_n)}(\phi_1), \ldots, E_n^{(h_n)}(\phi_\ell))$$
$$= n^{\ell/2} \mathrm{cum}(F_n^{(h_n)}(\phi_1), \ldots, F_n^{(h_n)}(\phi_\ell))$$
$$= \frac{1}{(2\pi)^\ell n^{\ell/2}} \sum_{t_1, \ldots, t_\ell} \sum_{k_1, \ldots, k_\ell} \hat\phi_1 \left( \frac{t_1}{n}, -k_1 \right) \cdots \hat\phi_\ell \left( \frac{t_\ell}{n}, -k_\ell \right)$$
$$\times \mathrm{cum}(X_{[t_1+1/2+k_1/2],n}^{(h_n)} X_{[t_1+1/2-k_1/2],n}^{(h_n)},$$
$$\ldots, X_{[t_\ell+1/2+k_\ell/2],n}^{(h_n)} X_{[t_\ell+1/2-k_\ell/2],n}^{(h_n)}).$$

We now use the representation $X_{t,n} = \sum_{j=-\infty}^{\infty} a_{t,n}(t-j) \varepsilon_j$ and obtain, with the product theorem for cumulants (cf. Brillinger (1981), Theorem 2.3.2) and the abbreviations $t_\nu^+ = t_\nu^+(t_\nu, k_\nu) = [t_\nu + 1/2 + k_\nu/2]$, $t_\nu^- = t_\nu^-(t_\nu, k_\nu) = [t_\nu + 1/2 - k_\nu/2]$, that this is equal to

$$\frac{1}{(2\pi)^\ell n^{\ell/2}}$$

$$\times \sum_{t_1,\ldots,t_\ell} \sum_{k_1,\ldots,k_\ell} \hat{\phi}_1\left(\frac{t_1}{n},-k_1\right)\cdots\hat{\phi}_\ell\left(\frac{t_\ell}{n},-k_\ell\right)$$

$$\times \sum_{i_1,\ldots,i_\ell,j_1,\ldots,j_\ell} \prod_{\nu=1}^{\ell}\left[h_n\left(\frac{t_\nu^+}{n}\right)h_n\left(\frac{t_\nu^-}{n}\right)a_{t_\nu^+,n}(t_\nu^+-i_\nu)\,a_{t_\nu^-,n}(t_\nu^--j_\nu)\right]$$

$$\times \sum_{\{P_1,\ldots,P_m\}\,\mathrm{i.p.}} \prod_{j=1}^{m}\mathrm{cum}(\varepsilon_s|s\in P_j),$$

where the last sum is over all indecomposable partitions (i.p.) $\{P_1,\ldots,P_m\}$ of the table

$$
\begin{array}{cc}
i_1 & j_1 \\
. & \\
. & \\
. & \\
i_\ell & j_\ell
\end{array}
\tag{63}
$$

with $|P_\nu| \geq 2$ (since $\mathbf{E}X(t)=0$). Using the upper bound $\sup_t |a_{t,n}(j)| \leq \frac{K}{\ell(j)}$ gives

$$
\begin{aligned}
&|\mathrm{cum}(E_n^{(h_n)}(\phi_1),\ldots,E_n^{(h_n)}(\phi_\ell))| \\
&\leq Kn^{-\ell/2}\sum_{t_1,\ldots,t_\ell}\sum_{k_1,\ldots,k_\ell}\left|\hat{\phi}_1\left(\frac{t_1}{n},-k_1\right)\cdots\hat{\phi}_\ell\left(\frac{t_\ell}{n},-k_\ell\right)\right| \\
&\qquad\times\prod_{\nu=1}^{\ell}\left[h_n\left(\frac{t_\nu^+}{n}\right)h_n\left(\frac{t_\nu^-}{n}\right)\right] \\
&\qquad\times\sum_{i_1,\ldots,i_\ell,j_1,\ldots,j_\ell}\prod_{\nu=1}^{\ell}\left[\frac{1}{\ell(t_\nu^+-i_\nu)}\frac{1}{\ell(t_\nu^--j_\nu)}\right] \\
&\qquad\times\sum_{\{P_1,\ldots,P_m\}\,\mathrm{i.p.}}\prod_{j=1}^{m}|\mathrm{cum}(\varepsilon_s|s\in P_j)|.
\end{aligned}
\tag{64}
$$

Concavity of $\log h_n(\cdot)$ implies that $h_n(\frac{t_\nu^+}{n})h_n(\frac{t_\nu^-}{n}) \leq h_n^2(\frac{[t_\nu+1/2]}{n})$. The Cauchy–Schwarz inequality (with squares of $\phi_1$ and $\phi_2$) now leads to the upper bound

$$
\begin{aligned}
Kn^{-\ell/2}&\left\{\sum_{t_1,\ldots,t_\ell}\sum_{k_1,\ldots,k_\ell}\left[\hat{\phi}_1\left(\frac{t_1}{n},-k_1\right)h_n^2\left(\frac{[t_1+1/2]}{n}\right)\right]^2\tilde{\phi}_3(k_3)\cdots\tilde{\phi}_\ell(k_\ell)\right. \\
&\qquad\times\sum_{i_1,\ldots,i_\ell,j_1,\ldots,j_\ell}\prod_{\nu=1}^{\ell}\left[\frac{1}{\ell(t_\nu^+-i_\nu)}\frac{1}{\ell(t_\nu^--j_\nu)}\right]
\end{aligned}
\tag{65}
$$

$$\times \sum_{\{P_1,\dots,P_m\}\,\mathrm{i.p.}} \prod_{j=1}^{m} |\mathrm{cum}(\varepsilon_s|s\in P_j)| \Bigg\}^{1/2} \Big\{ \text{ similar term } \Big\}^{1/2}$$

which, by using (51), is bounded by

$$Kn^{-\ell/2}\Bigg\{\sum_{t_1}\sum_{k_1,\dots,k_\ell}\left[\hat\phi_1\left(\frac{t_1}{n},-k_1\right)h_n^2\left(\frac{[t_1+1/2]}{n}\right)\right]^2\tilde\phi_3(k_3)\cdots\tilde\phi_\ell(k_\ell)$$

$$\times \sum_{i_1,\dots,i_\ell,j_1,\dots,j_\ell}\frac{1}{\ell(t_1^+-i_1)}\frac{1}{\ell(t_1^--j_1)}\prod_{\nu=2}^{\ell}\frac{1}{\ell(k_\nu-i_\nu+j_\nu)} \tag{66}$$

$$\times \sum_{\{P_1,\dots,P_m\}\,\mathrm{i.p.}}\prod_{j=1}^{m}|\mathrm{cum}(\varepsilon_s|s\in P_j)|\Bigg\}^{1/2}\Big\{ \text{ similar term } \Big\}^{1/2}.$$

Note that the term cum $(\varepsilon_s|s\in P_j)$ leads to the restriction that all $i_\nu, j_\nu \in P_j$ are equal. We now sum over the remaining indices from $k_2, i_1,\dots,i_\ell, j_1,\dots,j_\ell$, leading, due to the indecomposability of the partition and the fact that $1/\ell(j) \le K$, to the upper bound

$$Kn^{-\ell/2}\Bigg\{\sum_{t_1}\sum_{k_1,k_3,\dots,k_\ell}\left[\hat\phi_1\left(\frac{t_1}{n},-k_1\right)h_n^2\left(\frac{[t_1+1/2]}{n}\right)\right]^2\tilde\phi_3(k_3)\cdots\tilde\phi_\ell(k_\ell)\Bigg\}^{1/2}$$

$$\times\Bigg\{ \text{ similar term } \Bigg\}^{1/2}$$

$$\le Kn^{1-\ell/2}\rho_{2,n}^{(h_n)}(\phi_1)\,\rho_{2,n}^{(h_n)}(\phi_2)\prod_{j=3}^{\ell}\rho_\infty(\phi_j)$$

and therefore to the result.

In (ii), the generic constant $K$ needs to be independent of $\ell$. Again, we have (64) (with $K$ replaced by $K^\ell$). Remember that the term cum $(\varepsilon_s|s\in P_j)$ leads to the restriction that all $i_\nu, j_\nu \in P_j$ are equal. We denote this index by $i^{(j)}$ $(j=1,\dots,m)$.

We start by considering the case where $\ell$ is even. For each fixed partition $\{P_1,\dots,P_m\}$, we can renumber the indices $\{1,\dots,\ell\}$ in such a way that for every $j\in\{1,\dots,m\}$, there exists at least one even $\nu\in\{1,\dots,\ell\}$ and one odd $\nu\in\{1,\dots,\ell\}$ such that $i^{(j)}=i_\nu$ or $i^{(j)}=j_\nu$. This can be derived from the indecomposability of the partition and the fact that $|P_k|\ge 2$ for all $k$. The Cauchy–Schwarz inequality now yields as an upper bound for (64)

$$K^\ell\,n^{-\ell/2}\sum_{\{P_1,\dots,P_m\}\,\mathrm{i.p.}}\prod_{j=1}^{m}|\mathrm{cum}(\varepsilon_s|s\in P_j)|$$

$$\times \left\{ \sum_{t_1, \ldots, t_\ell} \sum_{k_1, \ldots, k_\ell} \prod_{\nu \text{ even}} \left[ \hat{\phi}_\nu \left( \frac{t_\nu}{n}, -k_j \right) h_n^2 \left( \frac{[t_\nu + 1/2]}{n} \right) \right]^2 \right.$$

$$\left. \times \sum_{i^{(1)}, \ldots, i^{(m)}} \prod_{\nu=1}^{\ell} \left[ \frac{1}{\ell(t_\nu^+ - i_\nu)} \frac{1}{\ell(t_\nu^- - j_\nu)} \right] \right\}^{1/2} \tag{67}$$

$$\times \left\{ \text{ the same term with } \ldots \prod_{j \text{ odd}} \ldots \right\}^{1/2},$$

where $i_\nu = i^{(j)}$ if $i_\nu \in P_j$ and $j_\nu = i^{(j)}$ if $j_\nu \in P_j$. By using relation (51), we have

$$\sum_{t_\nu, k_\nu} \frac{1}{\ell(t_\nu^+ - i_\nu)} \frac{1}{\ell(t_\nu^- - j_\nu)} \le K \sum_{k_\nu} \frac{1}{\ell(k_\nu - i_\nu + j_\nu)} \le K,$$

that is, the first bracket in (67) is bounded by

$$K^\ell \left\{ \sum_{t_\nu, k_\nu; \nu \text{ even}} \prod_{\nu \text{ even}} \left[ \hat{\phi}_\nu \left( \frac{t_\nu}{n}, -k_\nu \right) h_n^2 \left( \frac{[t_\nu + 1/2]}{n} \right) \right]^2 \right.$$

$$\left. \times \sum_{i^{(1)}, \ldots, i^{(m)}} \prod_{\nu \text{ even}} \left[ \frac{1}{\ell(t_\nu^+ - i_\nu)} \frac{1}{\ell(t_\nu^- - j_\nu)} \right] \right\}^{1/2}$$

$$\le K^\ell n^{\ell/4} \prod_{\nu \text{ even}} \rho_{2,n}^{(h_n)}(\phi_\nu).$$

The same applies for the second bracket in (67). Since the number of indecomposable partitions is bounded by $4^\ell (2\ell)!$, we obtain (62).

The case where $\ell$ is odd is a bit more involved. For each partition $\{P_1, \ldots, P_m\}$ with $m < \ell$, the result follows as in the case of $\ell$ even. For $m = \ell$, we can renumber the indices $\{1, \ldots, \ell\}$ such that each $P_\nu$ contains exactly one element of $\{i_\nu, j_\nu\}$ and $\{i_{\nu+1}, j_{\nu+1}\}$ (where $i_{\ell+1} = i_1, j_{\ell+1} = j_1$). For simplicity, we treat the case where $P_\nu = \{i_\nu, j_{\nu+1}\}$ ($\nu = 1, \ldots, \ell$) (the other cases follow analogously). We obtain, with the Cauchy–Schwarz inequality as an upper bound for each partition in (64),

$$K^\ell n^{-\ell/2} \sum_{t_\ell, k_\ell} \left\{ \sum_{t_1, \ldots, t_{\ell-1}} \sum_{k_1, \ldots, k_{\ell-1}} \prod_{j \text{ even}} \left[ \hat{\phi}_j \left( \frac{t_j}{n}, -k_j \right) h_n^2 \left( \frac{[t_j + 1/2]}{n} \right) \right]^2 \right.$$

$$\left. \times \sum_{i_1, \ldots, i_\ell} \prod_{\nu=1}^{\ell} \left[ \frac{1}{\ell(t_\nu^+ - i_\nu)} \frac{1}{\ell(t_\nu^- - i_{\nu-1})} \right] \right\}^{1/2}$$

$$\times \left\{ \text{the same term with} \ldots \prod_{j=1,3,\ldots,\ell-2} \ldots \right\}^{1/2},$$

where $i_0 = i_\ell$. By using relation (51), this is bounded by

$$K^\ell \left\{ \prod_{j=2}^{\ell-2} \rho_{2,n}^{(h_n)}(\phi_j) \right\} n^{-3/2} \sum_{t_\ell, k_\ell} \left| \hat{\phi}_\ell \left( \frac{t_\ell}{n}, -k_\ell \right) h_n^2 \left( \frac{[t_\ell + 1/2]}{n} \right) \right|$$

$$\times \left\{ \sum_{t_{\ell-1}, k_{\ell-1}} \frac{\hat{\phi}_{\ell-1}(\frac{t_{\ell-1}}{n}, -k_{\ell-1})^2 h_n^2(\frac{[t_{\ell-1}+1/2]}{n})}{\ell(t_\ell^- - t_{\ell-1}^+)} \right\}^{1/2}$$

$$\times \left\{ \sum_{t_1, k_1} \frac{\hat{\phi}_1(\frac{t_1}{n}, -k_1)^2 h_n^2(\frac{[t_1+1/2]}{n})}{\ell(t_\ell^+ - t_1^-)} \right\}^{1/2}.$$

The Cauchy–Schwarz inequality now yields, with (51),

$$K^\ell \left\{ \prod_{j=2}^{\ell-2} \rho_{2,n}^{(h_n)}(\phi_j) \right\} \rho_{2,n}^{(h_n)}(\phi_\ell)$$

$$\times n^{-1} \left\{ \sum_{t_1, k_1, t_{\ell-1}, k_{\ell-1}} \sum_{t_\ell, k_\ell} \frac{\hat{\phi}_{\ell-1}(\frac{t_{\ell-1}}{n}, -k_{\ell-1})^2 h_n^2(\frac{[t_{\ell-1}+1/2]}{n})}{\ell(t_\ell^- - t_{\ell-1}^+)} \right.$$

$$\left. \times \frac{\hat{\phi}_1(\frac{t_1}{n}, -k_1)^2 h_n^2(\frac{[t_1+1/2]}{n})}{\ell(t_\ell^+ - t_1^-)} \right\}^{1/2}$$

$$\leq K^\ell \left\{ \prod_{j=1}^{\ell} \rho_{2,n}^{(h_n)}(\phi_j) \right\},$$

which finally leads to (ii).

(iii) Follows from (i) since $\rho_{2,n}^{(h_n)}(\phi)^2 \leq K\rho_2(\phi)^2 + \frac{K}{n}\rho_\infty(\phi) \sup_j V(\hat{\phi}(\cdot, j))$. $\qquad \square$

**Proof of Theorem 2.5.** We have, for each $\phi_j$ (denoted by $\phi$ for simplicity) and $k \neq 0$,

$$\hat{\phi}(u, k) = \int_0^{2\pi} \frac{\exp(-\mathrm{i}k\lambda) - 1}{\mathrm{i}k} \phi_R(u, \mathrm{d}\lambda),$$

where $\phi_R(u, \mathrm{d}\lambda)$ is the signed measure corresponding to $\phi_R(u, \lambda) := \lim_{\mu \downarrow \lambda} \phi_R(u, \mu)$ (since $\phi$ is of bounded variation in $\lambda$, the limit exists; for the same reason, $\phi_R(u, \mathrm{d}\lambda)$ is a signed measure). This implies, for $k \neq 0$,

$$\sup_u |\hat{\phi}(u, k)| \leq \frac{K}{|k|} \sup_u V(\phi(u, \cdot)) = \frac{K}{|k|} \|\phi\|_{\infty, V} \quad \text{and}$$

$$V(\hat{\phi}(\cdot, k)) \leq \frac{K}{|k|} \|\phi\|_{V,V} \tag{68}$$

and, for $k = 0$,

$$\sup_u |\hat{\phi}(u, 0)| \leq 2\pi \|\phi\|_{\infty,\infty} \quad \text{and} \quad V(\hat{\phi}(\cdot, 0)) \leq 2\pi \|\phi\|_{V,\infty}.$$

Thus, $\rho_\infty(\phi)$ is not necessarily bounded and Theorem 5.3 cannot be applied. The trick now is to smooth $\phi(u, \lambda)$ in the $\lambda$-direction and to prove asymptotic normality instead for the resulting sequence of approximations: Let $k(x) := \frac{1}{\sqrt{2\pi}} \exp\{-\frac{1}{2}x^2\}$ be the Gaussian kernel, $k_b(x) := \frac{1}{b} k(\frac{x}{b})$ and

$$\phi_n(u, \lambda) = \int_{-\infty}^{\infty} k_b(\lambda - \mu)\, \phi(u, \mu)\, \mathrm{d}\mu$$

with $b = b_n \to 0$ as $n \to \infty$ (where $\phi(u, \mu) = 0$ for $|\mu| > \pi$). We have

$$\widehat{\phi_n}(u, k) = \hat{\phi}(u, k)\, \widehat{k_b}(k) \tag{69}$$

with $\widehat{k_b}(k) = \exp(-k^2 b^2/2)$. We obtain, from Lemma 5.7(i),

$$n \operatorname{var}[F_n(\phi_n) - F_n(\phi)] \leq \frac{K}{n} \sum_{t=1}^{n} \sum_{k=-\infty}^{\infty} \left( \widehat{\phi_n}\left(\frac{t}{n}, k\right) - \hat{\phi}\left(\frac{t}{n}, k\right) \right)^2$$
$$\leq K \sum_k \frac{[\exp(-k^2 b^2/2) - 1]^2}{k^2}. \tag{70}$$

Since $|1 - \exp(-k^2 b^2/2)| \leq \min\{1, \frac{k^2 b^2}{2}\}$, this is bounded by

$$K \sum_{|k| \leq 1/b} \frac{k^2 b^4}{4} + K \sum_{|k| > 1/b} \frac{1}{k^2} = O(b),$$

which implies that

$$\sqrt{n}(\{F_n(\phi_n) - \mathbf{E}F_n(\phi_n)\} - \{F_n(\phi) - \mathbf{E}F_n(\phi)\}) \xrightarrow{P} 0. \tag{71}$$

We now derive a CLT for $\sqrt{n}(F_n(\phi_{j,n}) - \mathbf{E}F_n(\phi_{j,n}))_{j=1,\dots,k}$ by applying Lemma 5.6(i) and Lemma 5.7(i). We obtain from (68) and (69) that

$$\sup_k |k|\, \tilde{\phi_n}(k) \leq K \|\phi\|_{\infty,V} \quad \text{and} \quad \sup_k |k| V(\hat{\phi_n}(\cdot, k)) \leq K \|\phi\|_{V,V}.$$

Furthermore, we have

$$\rho_\infty(\phi_n) \le 2\pi \|\phi\|_{\infty,\infty} + K\|\phi\|_{\infty,V} \sum_{k=1}^{\infty} \frac{1}{|k|} \exp(-k^2 b^2/2)$$

$$\le K(\|\phi\|_{\infty,\infty} + \log(b^{-1})\|\phi\|_{\infty,V})$$

and, with $h_n(\cdot) := I_{(0,1]}(\cdot)$,

$$\rho_{2,n}^{(h_n)}(\phi_n)^2 = \rho_{2,n}(\phi_n)^2 \le \rho_2(\phi_n)^2 + \frac{1}{n}\rho_\infty(\phi_n) \sup_k V(\hat{\phi_n}(\cdot,k))$$

$$\le \rho_2(\phi)^2 + O\left(\frac{\log(b^{-1})}{n}\right).$$

Therefore, the remainder term $R_n$ in Lemma 5.6(i) and the higher cumulants in Lemma 5.7(i) converge to zero if we choose, for example, $b = \frac{1}{n}$. Furthermore,

$$c_E^{(h_n)}(\phi_{j,n}, \phi_{k,n}) = c_E^{(I_{(0,1]})}(\phi_j, \phi_k) + O(b^{1/2}) \qquad (i, j = 1, \dots, k).$$

This follows easily by application of the Cauchy–Schwarz inequality and $\sup_u \int (\phi_{j,n}(u,\lambda) - \phi_j(u,\lambda))^2 \, d\lambda = O(b)$ (obtained with the Parseval formula as in (70)) and $\int (\phi_{j,n}(u,\lambda))^2 \, d\lambda \le K$. This gives the required CLT and, with (71), also the CLT for $\sqrt{n}(F_n(\phi_j) - \mathbf{E}F_n(\phi_j))_{j=1\dots,k}$. We obtain from (54) and (55) that

$$\sqrt{n}|\mathbf{E}F_n(\phi) - F(\phi)| \le K \frac{\|\phi\|_{\infty,\infty} + (\log n)\|\phi\|_{\infty,V} + \|\phi\|_{V,V} + \|\phi\|_{V,V}}{\sqrt{n}} = o(1),$$

which finally proves Theorem 2.5. $\qquad\qquad\square$

**Proof of Theorem 2.7 and Remark 2.8.** We obtain, from Lemma 5.7(ii) for the $\ell$th order cumulant in the case $\ell \ge 2$,

$$|\mathrm{cum}_\ell(E_n(\phi))| \le K^\ell (2\ell)! \rho_{2,n}(\phi)^\ell.$$

This implies, as on page 82 of Dahlhaus (1988), that

$$\mathbf{E}|(E_n(\phi))^\ell| \le (2K)^\ell (2\ell)! \rho_{2,n}(\phi)^\ell.$$

The result now follows in the same way as in the proof of Lemma 2.3 in Dahlhaus (1988).

We now prove the inequalities in Remark 2.8: We obtain, from the proof of Lemma 5.5 and an application of the Cauchy–Schwarz inequality,

$$\sqrt{n}|\mathbf{E}F_n(\phi) - F^+(\phi)| \le \sqrt{n}\rho_{2,n}(\phi)$$

$$\times \left(\frac{1}{n}\sum_{t=1}^{n}\left\{\sum_{|k|\le n}\left[\mathrm{cov}(X_{[t+1/2+k/2],n}^{(h_n)}, X_{[t+1/2-k/2],n}^{(h_n)}) - c\left(\frac{t}{n}, k\right)\right]^2\right.$$

$$+ \sum_{|k|>n} c\left(\frac{t}{n}, k\right)^2 \Bigg\}\Bigg)^{1/2}.$$

Application of Proposition 5.4 yields that the term in the bracket is of order $n^{-1/2}$, leading to (20). (21) and (22) follow by straightforward calculation, noting that Assumption 2.1 implies that $\sup_{u,\lambda} |f(u, \lambda)| \leq \infty$. (24) follows from an upper bound of (54) obtained by using (68). $\qquad\square$

**Proof of Theorem 2.9.** The proof uses a chaining technique as in Alexander (1984). Let

$$B_n = \left\{ \max_{t=1,\ldots,n} |X_{t,n}| \leq 2\log n \right\}. \tag{72}$$

Lemma 5.9 gives $\lim_{n\to\infty} P(B_n) = 1$. We will replace $\phi$ by

$$\phi_n^*(u, \lambda) = n \int_{u-1/n}^u \phi(v, \lambda)\, \mathrm{d}v \qquad (\text{with } \phi(v,\lambda) = 0 \text{ for } v < 0). \tag{73}$$

The reason for doing so is that otherwise we would need the exponential inequality (19) to hold with $\rho_2(\phi)$ instead of $\rho_{2,n}(\phi)$. Such an inequality does not hold. Instead, we exploit the following property of $\phi_n^*$:

$$\rho_{2,n}^2(\phi_n^*) = \frac{1}{n}\sum_{t=1}^n \int_{-\pi}^\pi \phi_n^*\left(\frac{t}{n}, \lambda\right)^2 \mathrm{d}\lambda = \frac{1}{n}\sum_{t=1}^n \int_{-\pi}^\pi \left(n\int_{(t-1)/n}^{t/n} \phi(u,\lambda)\,\mathrm{d}u\right)^2 \mathrm{d}\lambda$$

$$\leq \sum_{t=1}^n \int_{-\pi}^\pi \int_{(t-1)/n}^{t/n} \phi^2(u,\lambda)\,\mathrm{d}u\,\mathrm{d}\lambda = \rho_2^2(\phi). \tag{74}$$

Define

$$\tilde{E}_n^*(\phi) := \tilde{E}_n(\phi_n^*) = \sqrt{n}(F_n - \mathbf{E}F_n)(\phi_n^*). \tag{75}$$

Since the assertion and the proof of Theorem 2.7 are for $n$ fixed, we obtain from (19),

$$P(|\tilde{E}_n^*(\phi)| \geq \eta) \leq c_1 \exp\left(-c_2\sqrt{\frac{\eta}{\rho_2(\phi)}}\right). \tag{76}$$

On $B_n$, we have, by using Lemma 5.8 and Lemma 5.9, that

$$|\tilde{E}_n^*(\phi) - \tilde{E}_n(\phi)| \leq \sqrt{n}\,|F_n(\phi_n^*) - F_n(\phi)| + \sqrt{n}\,|\mathbf{E}F_n(\phi_n^*) - \mathbf{E}F_n(\phi)|$$

$$\leq 4\,K_1\left(\|\phi\|_{V,V}\,\frac{(\log n)^3}{n^{1/2}} + \|\phi\|_{V,\infty}\,\frac{(\log n)^2}{n^{1/2}}\right) + K_2\|\phi\|_{V,\infty}\frac{1}{n^{1/2}}$$

and therefore, with (26),

$$\sup_{\phi\in\Phi} |\tilde{E}_n^*(\phi) - \tilde{E}_n(\phi)| \leq \frac{\eta}{2}. \tag{77}$$

Thus,

$$P\left(\sup_{\phi \in \Phi} |\tilde{E}_n(\phi)| > \eta,\, B_n\right) \leq P\left(\sup_{\phi \in \Phi} |\tilde{E}_n^*(\phi)| > \frac{\eta}{2},\, B_n\right).$$

Let $\alpha := \tilde{H}_\Phi^{-1}(\frac{c_2}{4}\sqrt{\frac{\eta}{\tau_2}})$. We obtain, for any sequence $(\delta_j)_j$ with $\alpha = \delta_0 > \delta_1 > \cdots > 0$, where $\delta_{j+1} \leq \delta_j/2$ with $\eta_{j+1} := \frac{9}{c_2^2}\delta_{j+1}\tilde{H}_\Phi(\delta_{j+1})^2$,

$$\frac{\eta}{4} \geq \frac{18}{c_2^2} \int_0^\alpha \tilde{H}_\Phi(s)^2 \,\mathrm{d}s \geq \frac{18}{c_2^2}\sum_{j=0}^\infty (\delta_{j+1} - \delta_{j+2})\,\tilde{H}_\Phi(\delta_{j+1})^2 \geq \sum_{j=0}^\infty \eta_{j+1}. \tag{78}$$

For each number $\delta_j$, choose a finite subset $A_j$ corresponding to the definition of the covering numbers $N(\delta_j, \Phi, \rho_2)$. In other words, the set $A_j$ consists of the smallest possible number $N_j = N(\delta_j, \Phi, \rho_2)$ of midpoints of $\rho_2$-balls of radius $\delta_j$ such that the corresponding balls cover $\Phi$. Now, telescope

$$\tilde{E}_n^*(\phi) = \tilde{E}_n^*(\phi_0) + \sum_{j=0}^\infty \tilde{E}_n^*(\phi_{j+1} - \phi_j), \tag{79}$$

where the $\phi_j$ are the approximating functions to $\phi$ from $A_j$, that is, $\rho_2(\phi - \phi_j) < \delta_j$. The above equality holds on $B_n$ because Lemma 5.8 implies that

$$\sup_{\phi \in \Phi} |\tilde{E}_n^*(\phi - \phi_j)| \leq 5K_3\, n\, (\log n)^2 \sup_{\phi \in \Phi} \rho_2(\phi - \phi_j)$$

$$\leq 5K_3 n(\log n)^2\, \delta_j \to 0 \qquad \text{for all } n.$$

Thus,

$$P\left(\sup_{\phi \in \Phi} |\tilde{E}_n^*(\phi)| > \frac{\eta}{2},\, B_n\right)$$

$$\leq P\left(\sup_{\phi \in \Phi} |\tilde{E}_n^*(\phi_0)| > \frac{\eta}{4}\right) + \sum_{j=0}^\infty N_j N_{j+1} \sup_{\phi \in \Phi} P(|\tilde{E}_n^*(\phi_{j+1} - \phi_j)| > \eta_{j+1})$$

$$= I + II.$$

Hence, using exponential inequality (76), we have, by definition of $\alpha$, that

$$I \leq c_1 \exp\left\{ \tilde{H}_\Phi(\alpha) - \frac{c_2}{2}\sqrt{\frac{\eta}{\tau_2}} \right\} = c_1 \exp\left\{ -\frac{c_2}{4}\sqrt{\frac{\eta}{\tau_2}} \right\}. \tag{80}$$

In order to estimate $II$, we need a particular definition of the $\delta_j$. We set

$$\delta_{j+1} = \sup\left\{ x : x \leq \delta_j/2\,;\ \tilde{H}_\Phi(x) \geq \tilde{H}_\Phi(\delta_j) + \frac{1}{\sqrt{j+1}} \right\}.$$

Since $\sum_{\ell=1}^{j+1} \frac{1}{\sqrt{\ell}} \geq \int_0^{j+1} \frac{1}{\sqrt{x}} \, \mathrm{d}x = 2\sqrt{j+1} \geq 2\log(j+1)$, we obtain for term $II$

$$II \leq \sum_{j=0}^{\infty} c_1 \exp\left\{ 2\tilde{H}_\Phi(\delta_{j+1}) - c_2 \sqrt{\frac{9/c_2^2 \delta_{j+1} \, \tilde{H}_\Phi(\delta_{j+1})^2}{\delta_{j+1}}} \right\}$$

$$\leq \sum_{j=0}^{\infty} c_1 \exp\left\{ -\tilde{H}_\Phi(\delta_{j+1}) \right\} \leq \sum_{j=0}^{\infty} c_1 \exp\left\{ -\tilde{H}_\Phi(\alpha) - 2\log(j+1) \right\}$$

$$= c_1 \exp\left\{ -\frac{c_2}{4} \sqrt{\frac{\eta}{\tau_2}} \right\} \sum_{j=1}^{\infty} \frac{1}{j^2} \leq 2c_1 \exp\left\{ -\frac{c_2}{4} \sqrt{\frac{\eta}{\tau_2}} \right\}.$$

This implies the maximal inequality (28). To prove (29), we note that $|E_n(\phi) - \tilde{E}_n(\phi)| = \sqrt{n} \, |\mathbf{E}F_n(\phi) - F(\phi)|$, that is, instead of (77) on $B_n$, we obtain, with (81), (82), (21), (24) and (26),

$$\sup_{\phi \in \Phi} |\tilde{E}_n^*(\phi) - E_n(\phi)| \leq 13L \max\{\tau_{\infty, V}, \tau_{V, \infty}, \tau_{V, V}, \tau_{\infty, \infty}\} \frac{(\log n)^3}{\sqrt{n}} \leq \frac{\eta}{2}.$$

The rest of the proof is the same, that is, we also obtain (29). $\qquad \square$

**Proof of Theorem 2.11.** To prove weak convergence of $E_n$, we have to show weak convergence of the finite-dimensional distributions and asymptotic equicontinuity in probability of $E_n$ (cf. van der Vaart and Wellner (1996), Theorems 1.5.4 and 1.5.7). Convergence of the finite-dimensional distributions has been shown in Theorem 2.5. Asymptotic equicontinuity means that for every $\epsilon, \eta > 0$, there exists a $\tau_2 > 0$ such that

$$\liminf_n P\left( \sup_{\rho_2(\phi, \psi) < \tau_2} |E_n(\phi - \psi)| > \eta \right) < \epsilon.$$

In order to see this, we apply Theorem 2.9. For fixed $\eta > 0$, there exists a $\tau_2 > 0$ small enough such that (27) holds. To see this, notice that $\alpha \to 0$ as $\tau_2 \to 0$ and hence, using assumption (30), it follows that the integral on the right-hand side of (27) also tends to zero if $\tau_2 \to 0$. As $\eta > 0$ is fixed, (26) holds for $n$ large enough. Hence, we obtain, with $B_n$ from (72), for $\tau_2$ small enough,

$$\liminf_n P\left( \sup_{\rho_2(\phi, \psi) < \tau_2} |E_n(\phi - \psi)| > \eta \right)$$

$$\leq \liminf_n P\left( \sup_{\rho_2(\phi, \psi) < \tau_2} |E_n(\phi - \psi)| > \eta, \, B_n \right) + \lim_n P(B_n^c)$$

$$\leq 3c_1 \exp\left\{ -\frac{c_2}{4} \sqrt{\frac{\eta}{\tau_2}} \right\} < \epsilon. \qquad \square$$

**Proof of Theorem 2.12.** Let $\delta > 0$. The assumptions of Theorem 2.9 are fulfilled for $\eta = \delta\sqrt{n}$ and $n$ sufficiently large. For those $n$, we obtain, from Theorem 2.9,

$$P\left(\sup_{\phi \in \Phi_n} |F_n(\phi) - F(\phi)| > \delta\right) = P\left(\sup_{\phi \in \Phi_n} |E_n(\phi)| > \delta\sqrt{n}\right)$$

$$\leq 3c_1 \exp\left\{-\frac{c_2}{4}\sqrt{\frac{\delta\sqrt{n}}{\tau_2^{(n)}}}\right\} + P(B_n^c) \to 0. \qquad \square$$

**Lemma 5.8 (Properties of $F_n(\phi_n^*)$).** *Suppose Assumption 2.1 is fulfilled and $\phi_n^*(u,\lambda) = n\int_{u-1/n}^{u}\phi(v,\lambda)\,dv$ (with $\phi(v,\lambda) = 0$ for $v < 0$). We then have, with $X_{(n)} := \max_{t=1,\ldots,n}|X_{t,n}|$,*

$$|F_n(\phi) - F_n(\phi_n^*)| \leq K_1 X_{(n)}^2 \left(\|\phi\|_{V,V}\,\frac{\log n}{n} + \|\phi\|_{V,\infty}\,\frac{1}{n}\right), \tag{81}$$

$$|\mathbf{E}F_n(\phi) - \mathbf{E}F_n(\phi_n^*)| \leq K_2\|\phi\|_{V,\infty}\,\frac{1}{n}, \tag{82}$$

$$|F_n(\phi_n^*) - \mathbf{E}F_n(\phi_n^*)| \leq K_3(\sqrt{n}X_{(n)}^2 + 1)\rho_2(\phi). \tag{83}$$

**Proof.** We have, with (68),

$$|F_n(\phi) - F_n(\phi_n^*)| = \left|\frac{1}{n}\sum_{t=1}^{n}\int_{-\pi}^{\pi}\left(\phi\left(\frac{t}{n},\lambda\right) - \phi_n^*\left(\frac{t}{n},\lambda\right)\right)J_n\left(\frac{t}{n},\lambda\right)d\lambda\right|$$

$$\leq O(X_{(n)}^2)\sum_{t=1}^{n}\frac{1}{n}\sum_{k=-n}^{n}n\int_{(t-1)/n}^{t/n}\left|\widehat{\phi}\left(\frac{t}{n},-k\right) - \widehat{\phi}(u,-k)\right|du$$

$$\leq O(X_{(n)}^2)\frac{1}{n}\sum_{k=-n}^{n}V(\hat{\phi}(\cdot,k))$$

$$\leq K_1 X_{(n)}^2 \left(V^2(\phi)\frac{\log n}{n} + \sup_{\lambda}V(\phi(\cdot,\lambda))\frac{1}{n}\right).$$

Furthermore, we obtain, with Proposition 5.4 and (68),

$$|\mathbf{E}F_n(\phi) - \mathbf{E}F_n(\phi_n^*)| \leq \frac{1}{2\pi}\frac{1}{n}\sum_{t=1}^{n}\sum_{k=-n}^{n}n\int_{(t-1)/n}^{t/n}\left|\widehat{\phi}\left(\frac{t}{n},-k\right) - \widehat{\phi}(u,-k)\right|du$$

$$\times|\text{cov}(X_{[t+1/2+k/2],n}, X_{[t+1/2-k/2],n}|$$

$$\leq K_2\sup_{\lambda}V(\phi(\cdot,\lambda))\frac{1}{n}.$$

(83) has been proven in Lemma A.3 of Dahlhaus and Polonik (2006).                    $\square$

**Lemma 5.9.** *Suppose Assumption 2.1 is fulfilled with $\mathbf{E}|\varepsilon_t|^k \leq C_\varepsilon^k$ for all $k \in \mathbf{N}$. Let $X_{(n)} := \max_{t=1,\ldots,n} |X_{t,n}|$. We then have*

$$\mathbf{P}(X_{(n)} > 2\log n) = O(n^{-1}).$$

*That is, we have $\mathbf{P}(B_n^c) = O(n^{-1})$ for the set $B_n$ from (72) (used in Theorem 2.9).*

**Proof.** From (3), we have $\sup_{t,n} \sum_{j=-\infty}^{\infty} |a_{t,n}(j)| < m_0 < \infty$. The monotone convergence theorem and Jensen's inequality then imply

$$\mathbf{E}|X_{t,n}|^k \leq \mathbf{E}\left( \sum_{j=-\infty}^{\infty} |a_{t,n}(j)|\,|\varepsilon_{t-j}| \right)^k \leq m_0^k \mathbf{E}\left( \sum_{j=-\infty}^{\infty} \frac{|a_{t,n}(j)|}{m_0}\,|\varepsilon_{t-j}|^k \right) \leq m_0^k C_\varepsilon^k,$$

leading to

$$\mathbf{P}(X_{(n)} > 2\log n) \leq n \max_{t=1,\ldots,n} \mathbf{P}(X_{t,n} > 2\log n)$$

$$\leq n \frac{\mathbf{E}\,\mathrm{e}^{|X_{t,n}|}}{\mathrm{e}^{2\log n}} \leq \frac{1}{n} \sum_{k=0}^{\infty} \frac{m_0^k C_\varepsilon^k}{k!} \leq \frac{1}{n}\,\mathrm{e}^{m_0 C_\varepsilon} \to 0. \qquad \square$$

# Appendix: Proof of Proposition 2.4

We only give the proof for tvAR processes ($q = 0$). The extension to tvARMA processes is then straightforward. The proof is similar to that of Künsch (1995), who proved the existence of a solution of the form (2) under the assumption that the functions $\alpha_i(u)$ are continuous. Let

$$\boldsymbol{\alpha}(u) = \begin{pmatrix} -\alpha_1(u) & -\alpha_2(u) & \ldots & \ldots & -\alpha_p(u) \\ 1 & 0 & \ldots & \ldots & 0 \\ 0 & \ddots & \ddots & & \vdots \\ \vdots & \ddots & \ddots & \ddots & \vdots \\ 0 & \ldots & 0 & 1 & 0 \end{pmatrix}$$

and $\boldsymbol{\alpha}(u) = \boldsymbol{\alpha}(0)$ for $u < 0$. Since $\det(\lambda E_p - \boldsymbol{\alpha}(u)) = \lambda^p(\sum_{j=0}^{p} \alpha_j(u)\lambda^{-j})$, it follows that $\delta(\boldsymbol{\alpha}(u)) \leq \frac{1}{1+\delta}$ for all $u$ where $\delta(A) := \max\{|\lambda| : \lambda \text{ eigenvalue of } A\}$. Let

$$a_{t,n}(j) = \left( \prod_{\ell=0}^{j-1} \boldsymbol{\alpha}\left( \frac{t-\ell}{n} \right) \right)_{11} \sigma\left( \frac{t-j}{n} \right)$$

and

$$X_{t,n} = \sum_{j=0}^{\infty} a_{t,n}(j)\,\varepsilon_{t-j}.$$

It is easy to check that $X_{t,n}$ is a solution of (12) provided the coefficients are absolutely summable.

To prove this, we note (cf. Householder (1964), page 46) that for every $\varepsilon > 0$ and $u \in [0,1]$, there exists a matrix $M(u)$ with

$$\|\boldsymbol{\alpha}(u)\|_{M(u)} \leq \delta(\boldsymbol{\alpha}(u)) + \varepsilon, \tag{84}$$

where $\|A\|_M := \sup\{\|Ax\|_M : \|x\|_M = 1\}$ and $\|x\|_M = \|M^{-1}x\|_1 = \sum_{i=1}^{p} |(M^{-1}x)_i|$. Since the $\alpha_i(u)$ are functions of bounded variation (i.e., the difference of two monotonic functions), there exists for all $\varepsilon > 0$ a finite partition of intervals $I_1 \cup \cdots \cup I_m = [0,1]$ such that $|\alpha_i(u) - \alpha_i(v)| < \varepsilon$ for all $i$ whenever $u, v$ are in the same $I_k$. Let $M_k := M(u_k)$ for an arbitrary $u_k \in I_k$. Therefore, $m$ (and the partition) can be chosen such that

$$\|\boldsymbol{\alpha}(v)\|_{M_k} \leq \rho := \left(1 + \frac{\delta}{2}\right)^{-1} < 1 \qquad \text{for all } v \in I_k. \tag{85}$$

We now replace the first interval $I_1$ by $I_1 \cup (-\infty, 0)$ (remember that $\boldsymbol{\alpha}(u) = \boldsymbol{\alpha}(0)$ for $u < 0$). There exists a constant $c_0$ such that $\|B\|_1 := \sum_{i,j} |B_{i,j}| \leq c_0 \|B\|_{M_k}$ for all $k$. For $t$ and $n$ fixed, we now define $L_k := \{\ell \geq 0 : \frac{t-\ell}{n} \in I_k\}$ and $L_{k,j} := L_k \cap \{0, \ldots, j-1\}$. Then

$$
\begin{aligned}
|a_{t,n}(j)| &= \left| \left( \prod_{\ell=0}^{j-1} \boldsymbol{\alpha}\left(\frac{t-\ell}{n}\right) \right)_{11} \sigma\left(\frac{t-j}{n}\right) \right| \leq \left\| \prod_{\ell=0}^{j-1} \boldsymbol{\alpha}\left(\frac{t-\ell}{n}\right) \right\|_1 \sigma\left(\frac{t-j}{n}\right) \\
&\leq \prod_{k=1}^{m} \left\| \prod_{\ell \in L_{k,j}} \boldsymbol{\alpha}\left(\frac{t-\ell}{n}\right) \right\|_1 \sigma\left(\frac{t-j}{n}\right) \\
&\leq c_0^m \prod_{k=1}^{m} \left\| \prod_{\ell \in L_{k,j}} \boldsymbol{\alpha}\left(\frac{t-\ell}{n}\right) \right\|_{M_k} \sigma\left(\frac{t-j}{n}\right) \\
&\leq c_0^m \sup_u \sigma(u) \prod_{k=1}^{m} \rho^{|L_{k,j}|} = K\rho^j \qquad \text{(since } m \text{ is fixed)},
\end{aligned}
$$

that is, we have proven (3). Since $\|\boldsymbol{\alpha}(\frac{t-k}{n}) - \boldsymbol{\alpha}(\frac{t}{n})\|_1 = \sum_{i=1}^{p} |\alpha_i(\frac{t-k}{n}) - \alpha_i(\frac{t}{n})|$, we obtain, with similar arguments and

$$a(u,j) := (\boldsymbol{\alpha}(u)^j)_{11} \sigma(u),$$

$$
\begin{aligned}
&\sum_{t=1}^{n} \left| a_{t,n}(j) - a\left(\frac{t}{n}, j\right) \right| \\
&\leq \sum_{t=1}^{n} \sum_{k=1}^{j-1} \left\| \boldsymbol{\alpha}\left(\frac{t}{n}\right)^k \left( \boldsymbol{\alpha}\left(\frac{t-k}{n}\right) - \boldsymbol{\alpha}\left(\frac{t}{n}\right) \right) \prod_{\ell=k+1}^{j-1} \boldsymbol{\alpha}\left(\frac{t-\ell}{n}\right) \right\|_1 \sigma\left(\frac{t-j}{n}\right)
\end{aligned}
$$

$$+ \sum_{t=1}^{n} \left\| \boldsymbol{\alpha} \left( \frac{t}{n} \right)^j \right\|_1 \left| \sigma \left( \frac{t-j}{n} \right) - \sigma \left( \frac{t}{n} \right) \right|$$

$$\leq \sum_{t=1}^{n} \sum_{k=1}^{j-1} c_0 \rho^k \sum_{i=1}^{p} \left| \alpha_i \left( \frac{t-k}{n} \right) - \alpha_i \left( \frac{t}{n} \right) \right| c_0^m \rho^{j-1-k}$$

$$\leq K j^2 \rho^{j-1},$$

that is, (5). (4) and (6) follow similarly. □

# Acknowledgements

We are very grateful to two anonymous referees whose comments have led to substantial improvements of the paper. This work has been partially supported by the NSF Grant #0406431 and the Deutsche Forschungsgemeinschaft (DA 187/15-1).

# References

Alexander, K.S. (1984). Probability inqualities for empirical processes and the law of the iterated logarithm. *Ann. Probab.* **12** 1041–1067. Correction, *Ann. Probab.* **15** 428–430.

Brillinger, D.R. (1981). *Time Series: Data Analysis and Theory.* San Francisco: Holden Day. MR0595684

Dahlhaus, R. (1988). Empirical spectral processes and their applications to time series analysis. *Stochastic Process. Appl.* **30** 69–83. MR0968166

Dahlhaus, R. (1996). On the Kullback–Leibler information divergence of locally stationary processes. *Stochastic Process. Appl.* **62** 139–168. MR1388767

Dahlhaus, R. (1997). Fitting time series models to nonstationary processes. *Ann. Statist.* **25** 1–37. MR1429916

Dahlhaus, R. (2000). A likelihood approximation for locally stationary processes. *Ann. Statist.* **28** 1762–1794. MR1835040

Dahlhaus, R. and Neumann, M.H. (2001). Locally adaptive fitting of semiparametric models to nonstationary time series. *Stochastic Process. Appl.* **91** 277–308. MR1807680

Dahlhaus, R. and Polonik, W. (2002). Empirical spectral processes and nonparametric maximum likelihood estimation for time series. In *Empirical Process Techniques for Dependent Data* (H. Dehling, T. Mikosch and M. Sørensen, eds.) 275–298. Boston: Birkhäuser. MR1958786

Dahlhaus, R. and Polonik, W. (2006). Nonparametric quasi maximum likelihood estimation for Gaussian locally stationary processes. *Ann. Statist.* **34** 2790–2824. MR2329468

Davis, R.A., Lee, T. and Rodriguez-Yam, G. (2005). Structural break estimation for nonstationary time series models. *J. Amer. Statist. Assoc.* **101** 223–239. MR2268041

Fay, G. and Soulier, P. (2001). The periodogram of an i.i.d. sequence. *Stochastic Process. Appl.* **92** 315–343. MR1817591

Fryzlewicz, P., Sapatinas, T. and Subba Rao, S. (2006). A Haar–Fisz technique for locally stationary volatility estimation. *Biometrika* **93** 687–704. MR2261451

Householder, A.S. (1964). *The Theory of Matrices in Numerical Analysis.* New York: Blaisdell. MR0175290

Künsch, H.R. (1995). A note on causal solutions for locally stationary AR processes. Preprint, ETH Zürich.

Nason, G. P., von Sachs, R. and Kroisandt, G. (2000). Wavelet processes and adaptive estimation of evolutionary wavelet spectra. *J. Roy. Statist. Soc. Ser. B* **62** 271–292. MR1749539

Mikosch, T. and Norvaisa, R. (1997). Uniform convergence of the empirical spectral distribution function. *Stochastic Process. Appl.* **70** 85–114. MR1472960

Moulines, E., Priouret, P. and Roueff, F. (2005). On recursive estimation for locally stationary time varying autoregressive processes. *Ann. Statist.* **33** 2610–2654. MR2253097

Neumann, M.H. and von Sachs, R. (1997). Wavelet thresholding in anisotropic function classes and applications to adaptive estimation of evolutionary spectra. *Ann. Statist.* **25** 38–76. MR1429917

Priestley, M.B. (1965). Evolutionary spectra and non-stationary processes (with discussion). *J. Roy. Statist. Soc. Ser. B* **27** 204–237. MR0199886

Sakiyama, K. and Taniguchi, M. (2004). Discriminant analysis for locally stationary processes. *J. Multivariate Anal.* **90** 282–300. MR2081780

Van Bellegem, S. and Dahlhaus, R. (2006). Semiparametric estimation by model selection for locally stationary processes. *J. Roy. Statist. Soc. Ser. B* **68** 721–746. MR2301292

van der Vaart, A.W. and Wellner, J.A. (1996). *Weak Convergence and Empirical Processes.* New York: Springer. MR1385671

Whittle, P. (1953). Estimation and information in stationary time series. *Ark. Mat.* **2** 423–434. MR0060797